\documentclass[review]{elsarticle}
\usepackage{longtable}
\usepackage[utf8]{inputenc}
\usepackage{graphicx} 

\usepackage[normalem]{ulem} 
\usepackage{arydshln} 
\usepackage{amsfonts}                             
\usepackage{amssymb}                            
\usepackage{amsmath,amsthm} 
\usepackage[titletoc,toc,title]{appendix}
\usepackage{url}
\usepackage{pdflscape}
\usepackage{color}
\usepackage{graphicx,color}
\usepackage{nccmath}
\usepackage{subfig}
\usepackage{hyperref}
\usepackage{footnote}
\usepackage{algorithmicx}
\usepackage[ruled]{algorithm}
\usepackage{algpseudocode}
\usepackage{algc}
\usepackage{pdflscape}

\newtheorem{thm}{Theorem}
\newtheorem{prf}{Proof}

\begin{document}

\begin{frontmatter}

\title{A hybrid Lagrangian metaheuristic for the cross-docking flow shop scheduling problem}

\author[a]{Gabriela B. Fonseca}
\ead{gabrielabragafonseca@gmail.com}
\author[b]{Thiago H. Nogueira}
\ead{thiagoh.nogueira@ufv.br}
\author[a]{Mart\'in G\'omez Ravetti\corref{ca}}
\cortext[ca]{Corresponding author.}
\ead{martin.ravetti@dep.ufmg.br}

\address[a]{Department of Production Engineering, Universidade Federal de Minas Gerais, \\Av. Ant\^onio Carlos, 6627 ,CEP 31270-901 , Belo Horizonte, MG, Brazil.}
\address[b]{Department of Production Engineering, Universidade Federal de Vi\c{c}osa, \\Rodovia MG-230, CEP 38810-000, Rio Parana\'iba, MG, Brazil.}

\begin{abstract}
{
Cross-docking is a logistics strategy that minimizes the storage and picking functions of conventional warehouses.  The objective is to unload the cargo from inbound trucks and directly load it into outbound trucks, with little or no storage. The success of the strategy depends on an efficient transshipment operation. This work undertakes a study of truck scheduling in a parallel dock cross-docking center. The problem is first modeled as a two-machine flow shop scheduling problem with precedence constraints, with the objective of minimizing the makespan, and later we generalize it to the parallel-dock case. We propose a hybrid method based on a Lagrangian relaxation technique through the volume algorithm. Using information from the Lagrangian multipliers, constructive heuristics with local search procedures generates good feasible solutions. With a series of cuts, the methodology finds tight bounds for small and large instance sizes, outperforming current results.
}
\end{abstract}

\begin{keyword}
Logistics \sep Truck scheduling \sep Cross-docking \sep Lagrangian Relaxation.
\end{keyword}

\end{frontmatter}

\section{Introduction}

The current market environment, characterized by increasingly fierce competition, the globalization of the economy and an accelerated technological revolution has led companies to improve their production, logistics and distribution systems. Customers are demanding better and faster services, direct-to-customers strategies and many other new requirements, that expects more efficient and integrated logistic operations.  Cross-docking is one logistic technique that can reduce inventory costs while increasing the goods flow, improving the efficiency of the supply chain.

A Cross-docking Distribution Center (CDC) is a logistics technique widespread throughout the world. Several well-known companies such as retail chains (WalMart \cite{stalk1991}), mailing companies (UPS \cite{Forger1995}), automobile manufacturers (Toyota \cite{Witt1998}) and logistics providers (\cite{Gue1999}, \cite{Kim2008}) have gained considerable competitive advantages by using CDC.

The main idea behind cross-docking centers is to receive products from different suppliers or manufacturers and consolidate them to common final delivery destinations. In comparison to traditional warehouses, a CDC is managed with minimal handling and with little or no storage between unloading and loading of goods. This practice can serve different goals: the consolidation of shipments, a shorter delivery lead time, the reduction of costs, etc. Readers are referred to \cite{Ladier2016}, for a survey discussing industry practices and CDC problem characterization.

Cross-docks raise numerous optimization questions, either strategic, tactical or operational. In a Cross-docking Distribution Center (CDC), scheduling decisions are particularly relevant to ensure a rapid turnover and on-time deliveries. Due to its real-world importance, several truck-scheduling works and procedures have been introduced during recent years, treating specific cross-dock settings.

\begin{figure}[htbp]
    \begin{center}
        \includegraphics[scale=0.65]{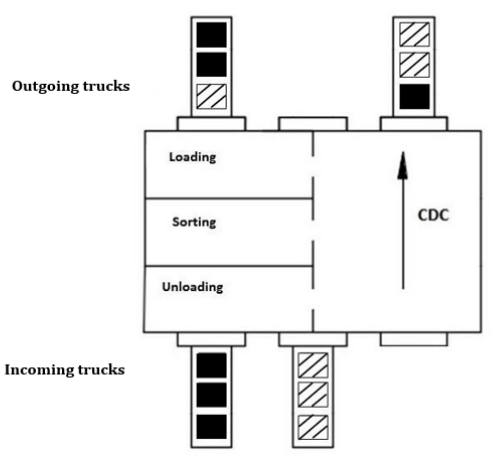}
        \caption{\label{Fig1}Schematic representation of a CDC.}
    \end{center}            
\end{figure}

A schematic representation of processes at a CDC is illustrated in Figure \ref{Fig1}. 
Firstly, incoming trucks arrive at the yard of the CDC, if the number of trucks is higher than the number of docks, some of them have to wait in a queue until further assignment. Secondly, after being docked, goods of the inbound trucks are unloaded, scanned, sorted, moved across the dock and loaded into outbound trucks for an immediate delivery elsewhere in the distribution chain. Once an outbound (inbound) truck is completely loaded (unloaded), the truck is removed from the dock, replaced by another truck and the course of action repeats.

In a conventional distribution center,  five operations are usually carried out when managing products: receiving, sorting, storing, picking and shipping operations. All of which can be done consecutively or in a particular order. As a result of applying a cross-docking system, a reduction of the cost of both storage and products picking can be achieved by synchronizing the inbound and outbound trucks flows. There are many benefits with the cross-docking: costs reduction (warehousing, inventory-holding, handling, and labor costs), shorter delivery lead times (from supplier to customer, increase of throughput), improved customer service and customer satisfaction, reduction of storage space, faster inventory turnover, fewer overstocks, reduced risk of loss and damage, etc. However, efficient transshipment processes and careful operations planning become indispensable within a CDC, where inbound and outbound flows need to be synchronized to keep the terminal storage as low as possible, and on-time deliveries are ensured. Many articles in the literature develop computerized scheduling procedures, which have achieved good results (See examples in Table \ref{Rev1}).

In this work, we consider a time-indexed integer programming formulation for a flow shop scheduling problem with cross-docking precedence constraints. We first consider the 2-dock case denoted as $F2|CD|C_{max}$, and strongly NP-hard (See \cite{ChenLee2009}), and later we deal with the parallel dock generalization, $F2(P)$ $|CD|C_{max}$. Computational experiments show that the model is efficient when solving small problems. 

A Hybrid Lagrangian Metaheuristic approach is also proposed and tested. The method can efficiently generate strong lower and upper bounds. Results are compared to the best heuristics in the literature introduced by Chen and Lee \cite{ChenLee2009} (based on Johnson's algorithm), with the CDH heuristic developed by Cota et al. \cite{CotaEtAll2016} for a parallel machines scenario, and to the best
heuristics proposed by Chen and Song \cite{ChenSong2009} for the multiple parallel processors (multiple docks) case.
Better solutions are consistently obtained for all instances tested. The article is organized as follows, an overview of related works is presented in Section 2. The mathematical model and the Hybrid Lagrangian Metaheuristic Framework are described in Section 3 and 4, respectively. Section 5 generalizes the algorithm for the parallel machine case. Computational experiments are reported in Section 6. Finally, discussions and conclusions are drawn in Section 7. 


In this manuscript, we propose a Hybrid Lagrangian Metaheuristic Framework for the cross-docking flow shop scheduling problem. The work has a similar line of thought than Boschetti and Maniezzo \cite{Boschetti2009}. However, it emphasizes the use of the information of the Lagrangian Multipliers as discussed in Pirkwieser et al. \cite{PirkwieserEtAl2007}, and Paula et al. \cite{PaulaEtAl2010}.

For the 2-dock case, although the subproblems of the Lagrangian relaxation are polynomial, through a series of cuts on the makespan value, we improve the linear lower bounds. The proposed algorithm uses the information obtained from the Lagrangian multipliers to construct feasible solutions through a NEH based heuristic (\cite{Nawaz1983}), and performs search procedure through a Local Search framework. Furthermore, the proposed hybrid algorithm proves optimality in several instances.

We generalize our approach to the parallel-dock case. In this scenario, the subproblems are NP-hard, with a focus on improving the feasible solutions we solve the subproblem's linear relaxations, to feed the Lagrangian Metaheuristic.  Results show that the technique outperforms current heuristic methods.

\section{Literature Review}

The increasing use of cross-docking techniques has motivated several authors to investigate new ways to improve the design and tactical operations of CDCs. In the literature, several decision problems are studied concerning strategical, tactical, and operational decisions. 

Boysen and Fliedner \cite{BoysenFliedner2010}, and Belle et al. \cite{Belle2012} focus on the review of the existing literature of cross-docking problems. According to \cite{Boysen2010}, decision problems in a cross-docking center can be allocated according to the following classification, ordered from strategic to operational levels: location of cross-docking centers, layout, vehicle routing, dock's assignments of destinations, truck scheduling and resource scheduling inside the center. Here we focus on truck scheduling problems, where the aim is deciding where and when should the trucks be processed.

Examples of strategic decisions in a CDC can be found in many works. Location of a cross-docking is analyzed in Campbell \cite{Campbell1994}, Klose and Drexl \cite{Klose2005}, Chen et al. \cite{Chenetal2006}, and Chen \cite{Chen2007}. The layout of cross-docking centers is studied in Gue \cite{Gue1999}, Bartholdi and Gue \cite{BartholdiGue2002} and Vis and Roodbergen \cite{Vis2008}. 

Some other works deal with the vehicle routing as Oh et al. \cite{Oh2006} and Wen et al. \cite{Wen2007}. Regarding operational decisions, some works consider the dock assignment problem: Belle et al. \cite{Belle2012}, Boysen and Fliedner \cite{BoysenFliedner2010}, Tsui and Chang \cite{Tsui1990} \cite{Tsui1992}, Bozer and Carlo \cite{Carlo2008} and Gue \cite{Gue1999}.

In a CDC, scheduling decisions are particularly important to ensure a rapid turnover and on-time deliveries. Due to its real-world importance, several truck-scheduling works and procedures have been introduced during recent years, treating specific cross-dock settings. For that reason, we highlight in Table \ref{Rev1}, some outstanding works dealing with truck scheduling problems in a CDC. As proposed by \cite{Boysen2010}, a classification scheme for deterministic truck scheduling problems is presented for each work, as well as, their contribution.

\newpage
\begin{table}[htbp]
  \caption{Previous specific research works for the truck scheduling problems. Based on \cite{BoysenFliedner2010}.}  
	 \begin{center}
	 \scalebox{0.75}{
	\begin{tabular}{c c c c} 
		\hline
		\multicolumn{1}{c}{\textbf{Publication}} & \multicolumn{1}{c}{\textbf{Notation}} & \multicolumn{1}{c}{\textbf{Complexity}} & \multicolumn{1}{c}{\textbf{Contribution}} \\
		\hline
		Miao et al. \cite{Miao2009}& [$M|limit, t_{io}|\star$] &NP-hard&MM, HM, P\\
		Chen and Lee \cite{ChenLee2009}& [$E2|t_{j}=0|C_{max}$] &NP-hard&B, ES, P\\
		Chen and Lee \cite{ChenLee2009}& [$E2|t_{j}=0|C_{max}$] &NP-hard&P\\
		Chen and Song \cite{ChenSong2009}& [$E|t_{j}=0|C_{max}$] &NP-hard&MM, HS, B\\
		Boysen \cite{Boysen2010} & [$E|p_{j}=p, no-wait, t_{j}=0|\sum T_{o}$] & Open & MM, HM, ES \\
		Boysen and Fliedner \cite{BoysenFliedner2010}& [$E|t_{io}, f_{ix}|\sum W_{s}U_{s}$] &NP-hard&MM, P\\
		Boysen and Fliedner \cite{BoysenFliedner2010}&[$E|t_{i}= 0, f_{ix}$ $|$ $\sum W_{s}U_{s}$] &NP-hard&P\\
		Boysen et al. \cite{BoysenEtAll2010}& [$E2| p_{j}=p, change|C_{max}$] &NP-hard&MM, HS, HI, ES, P\\
		McWilliams \cite{Mcwilliams2010}& [$E|no-wait|*$] &Open&HM\\			
		Vahdani and Zandieh \cite{Vahdani2010}& [$E2| change |C_{max}$] &NP-hard& MM, HM\\	
		Melo and Araujo \cite{AraujoeMelo2010}& [$E|t_{j}=0|C_{max}$] &NP-hard& HM\\		
		Arabani et al. \cite{Arabani2011}& [$E2|change|C_{max}$] &Open& HM\\
		Larbi et al. \cite{Larbi2011}& [$E2|pmtn|*$] &NP-hard&MM, ES\\			
		Alpan et al. \cite{Alpan2011}& [$E|pmtn|*$] &Open&MM, ES\\						
		Lima \cite{Lima2014}& [$E2|t_{j}=0|\sum C^{2}_{j}$] &NP-hard& MM, HS, HM\\
		Cota et al. \cite{CotaEtAll2016}&[$E|t_{j}=0|C_{max}$]&NP-hard& MM, HS\\
	\hline
\end{tabular}
}

\footnotesize{The notation used in column ``Contribution'' is stated as: MM (mathematical model), HI (heuristic improvement procedure), HM (meta-heuristic), B (bound computation), HS (start heuristic for initial solution), ES (exact solution procedure), P (properties [e.g., complexity] of problem). In the Boysen's \cite{Boysen2010} tuple notation, column ``Notation'', the fields are door environment, operational characteristics and the objective, respectively.
}
\end{center}
  \label{Rev1}
\end{table}

Chen and Lee \cite{ChenLee2009} study a two-machine cross-docking flow shop scheduling problem, in which a job at the second machine can be processed only after finalizing some jobs at the first one, with the objective of minimizing the makespan. The authors show the problem is strongly NP-hard. They develop a polynomial approximation algorithm with an error-bound analysis and a branch-and-bound algorithm. Computational results show that the branch-and-bound algorithm can optimally solve problems with up to 60 jobs in a reasonable time.

Chen and Song \cite{ChenSong2009} extended the Chen and Lee \cite{ChenLee2009} problem to the two-stage hybrid cross-docking scheduling problem, by considering multiple parallel processors (multiple docks) per stage (inbound and outbound), allowing simultaneous loading and unloading operations. They propose a mixed integer programming model and four constructive heuristics, based on Johnson's rule, to investigate the performance for moderate and large scale instances.

Another studies considering multiple dock environment are proposed by Alpan et al. \cite{Alpan2008}, but with a different aim, a dynamic program based model is presented to find the schedule of transshipment operations minimizing the operations total cost, solving only small size instances. Larbi et al. \cite{Larbi2009} presented some heuristics to obtain near optimal solutions for large instances considering the same scenario than~\cite{Alpan2008}.

Cota {\it et al.} \cite{CotaEtAll2016} deals with the operational decision problem of scheduling the trucks on multiple inbound and outbound docks, that is, 
they consider the same setting than Chen and Song \cite{ChenSong2009}. They propose a time-indexed mixed integer linear programming formulation and 
a constructive polynomial heuristic. Results are compared to the best heuristics proposed in \cite{ChenSong2009,CotaEtAll2016} on medium and large instance, and better solutions were consistently obtained for instances with greater number of machines and jobs.

The hybrid metaheuristics arise as an option to achieve solutions within a reasonable computational time in very difficult combinatorial optimization problems. The hybrid method consists of a cooperative combination of methods, exact and/or approximated, and aims to absorb to the limit the potentialities of all the approaches (see surveys Alba \cite{Alba2005}, Puchinger and Raidl \cite{Puchinger2005} and Blum et al. \cite{BlumEtAl2011}). 

Paula {\it et al.} \cite{PaulaEtAl2010} propose a variant of the Lagrangian Relaxation to obtain good bounds to the problem of parallel machines scheduling with sequence-dependent setup times. In this effort, the Lagrangian Relaxation is improved by the use of a metaheuristic as an internal procedure of the Lagrangian Relaxation. The heuristic is used to generate feasible solutions during the execution of the non-delayed relax-and-cut algorithm. 

Pirkwieser {\it et al.} \cite{PirkwieserEtAl2007} present a similar work, using the Lagrangian Relaxation applied to the
Knapsack Constrained Maximum Spanning Tree problem, in which the Lagrangian Multipliers
influence the procedures executed by the Genetic Algorithm.

Boschetti and Maniezzo \cite{Boschetti2009} propose a Lagrangian Metaheuristic procedure. Their approach consists in using metaheuristic techniques to obtain feasible solutions using the information of Lagrangian Multipliers. The solutions obtained tend to be good bounds. We remark that, in these efforts, the potential of the information of the Lagrangian Multipliers is still not completely used. Part of our goals is to improve the use of this information. Readers are referred to \cite{BlumEtAl2011} for a survey on hybrid metaheuristics. 

As discussed here, a significant number of CDC articles are dedicated to the study of strategical and tactical aspects. However, not always aligned with industry practice. More recently, Ladier and Alpan~\cite{Alpan2011} discuss the relationship between the industry needs and academic research topics, pointing out the gap between both worlds. They highlight some gaps between theory and practice which would need to be filled in order to get closer to real-life constraints, for example: researchers should aim towards removing simplifying assumptions which are difficult to justify in an industrial context, or, it would be necessary to take deadlines into account for outbound truck departures.

In this work, we consider the problem of scheduling the trucks as described in Chen and Lee's approach \cite{ChenLee2009}. As Chen and Song  \cite{ChenSong2009} and Cota {\it et al.} \cite{CotaEtAll2016} generalize their problem, we compare our results to the best algorithms in all these works considering a 2-dock version and a parallel-dock case.

\section{Mathematical Model}\label{model}

In this section, we present a mathematical formulation for a 2-dock version based on scheduling problem initially proposed by Chen and Lee \cite{ChenLee2009}. 
The integer programming model considered in this work adopts a time-indexed formulation proposed by Lima \cite{Lima2014} but adopting the makespan as  objective 
function, as proposed by Chen and Lee \cite{ChenLee2009}. 

Considering a 2-dock cross-docking center where $n$ loaded trucks arrive with products demanded by one or more customers. In the center, each truck must unload its cargo and load it into one or more outbound trucks, responsible for delivery to specific destinations in the supply chain. Each departing vehicle may leave the center only after the cargo is fully loaded, and it can begin the loading process only after all needed products were available in the center. Here we consider the existence of two docks in the center, one dock dedicated to unload the trucks and the other one for loading purposes, machine 1 (M1) and machine 2 (M2), respectively. 
The problem corresponds to define the truck sequence in the inbound and outbound dock, minimizing the completion time of the last job processed by machine 2.

The following parameters, sets and variables are considered:\\

\textbf{Input Parameters} 

\begin{itemize}
\item $n$: number of jobs to be processed on M1.
\item $m$: number of jobs to be processed on M2.
\item $p_{ij}$: processing time of job $j$ on machine $i$.
\item $T_f$: size of time horizon. A first estimate for the time horizon is the sum of the processing times of all jobs.
\item $T_0 \geq 0$: beginning of time horizon for the jobs ${j \in J^2}$, i.e., minimum starting date for processing of the output trucks, initially set as $T_0$ = $min_{j \in J^2}$ $\{\sum_{i \in J^1,i \in S_j} p_{1i}\}$, then $T$ = \{$T_0,...,T_f$\} $\forall j \in J^2$.
\end{itemize}

\textbf{Sets}\\

The set of periods is defined as $T$ = \{$0,...,T_f$\}. To represent arrivals and truck departures in cross-docking center, two sets of jobs are created:
\begin{itemize}

\item  $J^1$ = \{$j_1^1, j_2^1,...,j_n^1$\} , set of $n$ jobs (inbound trucks), which must be processed on M1.

\item  $J^2$ = \{$j_1^2, j_2^2,...,j_m^2$\}, set of $m$ jobs (outbound trucks), which must be processed on M2.

\item  $S_j$, set of precedent jobs. For each job $j_j^2 \in J^2$, there is a corresponding subset of $J^1$ that must be completed before beginning its process.  It is considered that each element job ${j \in J^2}$ has at least one job $S_j$ as precedent.\\
\end{itemize}

\textbf{Decision variables} 

\begin{itemize}
\item $x_{jt}$ $(\forall j \in J^1, \forall t \in T$), binary variable, $x_{jt}$ is equal to 1 if job $j$ starts its process at time $t$ and equal to $0$ otherwise.
\item $y_{jt}$ $(\forall j \in J^2, \forall t \in T$), binary variable, $y_{jt}$ is equal to 1 if job $j$ starts its process at time $t$ and equal to $0$ otherwise.
\item $C_{max}$: makespan. Maximum completion time in M2.
\end{itemize}

\newpage
The complete mathematical model (F) is presented as follows:

\begin{fleqn}
\begin{equation}
(F)\hspace{2 cm} {Minimize}   \hspace{0.2cm} {C_{max} \label{foeq}} 
\end{equation}
\end{fleqn}

subject to\\
\begin{eqnarray}
\sum_{t=0}^{T_f-p_{1j}} x_{jt} = 1, & & \forall j \in J^1, \label{eq01} \\
\sum_{t=T_0}^{T_f-p_{2j}} y_{jt} = 1, & & \forall j \in J^2, \label{eq02} \\
\sum_{j \in J^1} \sum_{s=max{(0 ; t-p_{1j}+1)}}^t x_{js} \leq 1, & & \forall t \in T, \label{eq03}\\
\sum_{j \in J^2} \sum_{s=max{(T_0 ; t-p_{2j}+1)}}^t y_{js} \leq 1, & & \forall t \in T, \label{eq04}\\
\sum_{t=T_0}^{T_f-p_{2j}} ty_{jt} - \sum_{t=0}^{T_f-p_{1k}} (t+p_{1k})x_{kt} \geq 0, & &  \forall j \in J^2, \forall k \in S_j, \label{eq05}\\
C_{max} \geq p_{2j} + \sum_{t=T_0}^{T_f-p_{2j}}ty_{jt}, & & \forall j \in J^2, \label{eq06}\\
x_{jt} \in \{0,1\} , & &  \forall j \in J^1, \forall t \in T, \label{eq07}\\
y_{jt} \in \{0,1\} , & &  \forall j \in J^2, \forall t \in T, \label{eq08}\\
C_{max} \geq 0.\label{eq09}
\end{eqnarray}

The objective function (\ref{foeq}) minimizes the makespan. The set of constraints (\ref{eq01}) ensures that each job $j_j^1$ $\in$ $J^1$ should start its processing in one and only one period within the time horizon. The set (\ref{eq02}) works with the same reasoning on M2. Constraints (\ref{eq03}) ensure that a job $j_j^1$ $\in$ $J^1$ does not start its processing while another job is being processed in the same machine M1. The constraints set (\ref{eq04}) works, in the same way, but applied to jobs $j_j^2$ $\in$ $J^2$. The set (\ref{eq05}) controls the precedence relationships. For every existing precedence relation, the start date of the job $j_j^2$ $\in$ $J^2$ should be greater than or equal to the completion time of its precedent $j_k^1 \in S_j$. The set of constraints (\ref{eq06}) indicates that the variable $C_{max}$ should be the maximum completion time of jobs in $J^2$. Finally, sets (\ref{eq07}) to (\ref{eq09}) define the variable's domain.

\subsection{Lagrangian Relaxation}
\label{HgR}


Let  L($\lambda$) be the relaxation of formulation F, when constraints \ref{eq05} are dualized and its violation penalized in the objective function. For each constraint is associated one Lagrangian multiplier $\lambda$ representing the weight given to the violation. The set of Lagrangian multipliers will be represented by $\lambda_{jk}$, with $j \in J^2$ and 
$k \in S_j$. Notice that, these constraints are the cross-docking precedence relations, coupling decisions in M1 with decisions in M2.
Thus, the subproblem L($\lambda$) is defined as:

\begin{fleqn}
\begin{equation}
 L(\lambda) = Minimize \hspace{0.2cm} {C_{max}} + \sum_{j \in J^2} \sum_{k \in S_j} \lambda_{jk} ( \sum_{t=0}^{T_f-p_{1k}} (t+p_{1k})x_{kt} - \sum_{t=T_0}^{T_f-p_{2j}} ty_{jt} ) \label{eqsublag}
\end{equation}
\end{fleqn}

s.t. (\ref{eq01}) -- (\ref{eq04}), (\ref{eq06}) -- (\ref{eq09}), $\lambda_{jk}  \geq 0$.



Now we are able to uncouple the problem into two new scheduling problems, one for each machine. Subproblem $X$ considers M1 while subproblem $Y$ considers M2. The value of the lower bound is obtained by adding the subproblems objective functions L$(\lambda)_x$ and L$(\lambda)_y$.

\textbf{Subproblem $X$}

Isolating the terms of the objective function that contain variables $x$, we reach to:

\begin{equation}
 L(\lambda)_x = Minimize \sum_{j \in J^2} \sum_{k \in S_j} \lambda_{jk} \sum_{t=0}^{T_f -p_{1k}} (t+p_{1k})x_{kt}  \label{eqsublagx}
\end{equation}
s.t.  (\ref{eq01}), (\ref{eq03}), (\ref{eq07}) and $\lambda_{jk}  \geq 0$.


Rewriting the above objective function from the perspective of jobs in \textit{$J^1$}, we have: 

\begin{equation}
L(\lambda)_x  = Minimize \sum_{j \in J^1} \sum_{t=0}^{T_f -p_{1j}} (t+p_{1j})x_{jt} w_j^1 \label{eqsublagx2}
\end{equation}

\begin{flushleft}
where $\quad w_j^1 = \sum_{i \in J^2} \lambda_{ij}$,  where $ \lambda_{ij} = 0$ if $\textit{j} \notin S_i, \quad$.
\end{flushleft}

%

The problem mentioned above is known as The Total Weighted Completion Time (\cite{Pinedo2008}), denoted by  $1 ||\sum C_j W_j$. This problem can be solved by the WSPT rule (Weighted Shortest Processing Time First)  proposed by Smith in 1956 \cite{Smith1956}. According to this rule, the optimal solution is obtained by ordering the jobs in descending order of \mbox{$w_j^1$ / $p_{1j}$}, and can be obtained in O($n\log$($n$)).\\

\textbf{Subproblem $Y$}

With the terms associated to variables  $y$, subproblem $Y$ can be defined as:\\

\begin{equation}
 L(\lambda)_y = Minimize \hspace{0.2cm} {C_{max}} - \sum_{j \in J^2} \sum_{k \in S_j} \lambda_{jk} \sum_{t=T_0}^{T_f-p_{2j}} ty_{jt} \label{eqsublagy}
\end{equation}
subject to\\
\begin{eqnarray}
\sum_{t=T_0}^{T_f-p_{2j}} y_{jt} = 1, & & \forall j \in J^2, \label{eqy02} \\
\sum_{j \in J^2} \sum_{s=max{(T_0 ; t-p_{2j}+1)}}^t y_{js} \leq 1, & & \forall t \in T, \label{eqy04}\\
C_{max} \geq p_{2j} + \sum_{t=T_0}^{T_f-p_{2j}}ty_{jt}, & & \forall j \in J^2, \label{eqy06}\\
y_{jt} \in {0,1}, & &  \forall j \in J^2, \forall t \in T, \label{eqy08}\\
\lambda_{jk}  \geq 0, & & \forall j \in J^2, k \in S_j, \label{eqy09}\\
C_{max} \geq 0. \label{eqy10}
\end{eqnarray}

Defining $w_j^2 = - \sum_{k \in S_j} \lambda_{jk}$ $\forall$$j$ $\in$ $J^2$, we can rewrite the problem as follows:

\begin{equation}
L(\lambda)_y = Minimize \hspace{0.2cm} {C_{max}} + \sum_{j \in J^2} \sum_{t=T_0}^{T_f-p_{2j}} t.y_{jt} w_j^2  \label{eqsublagy2}
\end{equation}
s.t.
(\ref{eqy02}), (\ref{eqy04}), (\ref{eqy06}), (\ref{eqy08}) e (\ref{eqy10}).

The first term of the objective function~(\ref{eqsublagy2}) is the makespan. The second term represents the sum of the weighted starting times of the jobs on M2, also known as The Total Weighted Starting Time, denoted $\sum I_j W_j$. Thus, the Y subproblem can be defined as $1||{C_{max}} +\sum I_j W_j$. To solve this subproblem, we proposed a rule named WSPT-TRD. 

The WSPT-TRD rule works as follows. First, the WSPT rule is used to generate a sequence of jobs on M2. The second step is to define the correct allocation of jobs in the time horizon. For that, we evaluate the increase of the makespan when the weights associated with the jobs are negative. When the weights are negative, it creates a trade-off situation, as explained next.

When considering a general, weights $w_j^2$ may have positive or negative values. Let us first consider schedule decisions for jobs with $w_j^2 \geq 0$, both $C_{max}$ and $\sum I_j W_j$ have similar behavior since the objective function is minimized by allocating the jobs at the beginning of the time horizon. In our particular case, as $\lambda_{jk}  \geq 0$, we do not have the case of $w_j^2 > 0$; however, our algorithm here proposed considers this possibility.

Considering jobs with negative weights ($w_j^2 < 0$), the objective function terms lead to different solutions. On the one hand, by considering the makespan $C_{max}$ one tends to allocate the jobs at the beginning of time horizon, as the weights do not influence its final value. On the other hand, when analyzing the $\sum I_j W_j$ one tends to allocate the jobs as late as possible, taking advantage of the negative weights. 

We notice that if we delay all jobs with negative weights one-time unit, the $C_{max}$ increases one unit.  Thus, our algorithm computes the sum of the negative weights, if the sum of the negative weights is less than $-1$, the displacement of the jobs towards the end of the time horizon compensates the increase of $C_{max}$, defining in this way the best allocation of jobs.

\subsection*{Pseudocode of algorithm WSPT-TRD:}

\begin{description}
\item{\textit{Step 1:}} Compute the weights $w_j^2$ on the machine 2 ($w_j^2 = - \sum_{k \in S_j} \lambda_{jk} \hspace{0.1cm}\forall j \in J^2$).
\item{\textit{Step 2:}} Sort the jobs in $J^2$ in decreasing order of $w_j^2$ / $p_{2j}$.

\item{\textit{Step 3:}} Compute the sum of the negative weights as $\delta$ = $\sum_{j \in J^2}w_j^2 \hspace{0.1cm} \forall  j \in J^2 \wedge w_j^2<0$.


\item{\textit{Step 4:}} If $w_j^2 = 0$ $\forall j \in J^2$ or $\delta$ $\geq -1$, allocate the jobs from $t=T_0$ by sequence generated by WSPT rule. 
\item{\textit{Step 5:}} If there exists $w_j^2 < 0$ and $\delta$ $\leq -1$, allocate the jobs from t = $T_f$ by the reverse sequence generated by WSPT rule.
\end{description}

\begin{thm}
Considering the problem $1||{C_{max}} +\sum I_j W_j$, the algorithm WSPT-TRD obtains an optimal solution.
\end{thm}
\begin{prf}
The proof is provided in Appendix \ref{proof1}.
\end{prf}

\subsection{Lower Bounds}

Besides the lower bound obtained from the Lagrangian relaxation, we consider here two other lower bounds. The first one, called $LB_1$ is defined as proposed by Chen and Lee \cite{ChenLee2009}, as follows:

\begin{equation}
LB_1 = \sum_{i \in J^2} p_{2i} + min_{j \in J^2} \{\sum_{i \in J^1,i \in S_j} p_{1i}\}.
\end{equation}

$LB_1$ basically computes the sum of processing time on $M2$, plus the minimum amount of time necessary to begin the process of the first job in $M2$, that is, the minimum set of precedent jobs in $M1$.

The second lower bound, $LB_2$, proposed in this article, considers a complementary condition:

\begin{equation}
LB_2 = \sum_{i \in J^1} p_{1i} + min_{j \in J^1} \{\sum_{i \in J^2,i \in P_j} p_{2i}\}.
\end{equation}

Where, $P_j$ represents the successors subsets jobs $j \in J^1$ corresponding to job $i \in J^2$, i.e., successors are all jobs that must be processed after completion for a given job in the first stage.

In this case, $LB_2$ computes the sum of processing times in $M1$ plus the minimum set of common successors of $j \in J^1$. Finally, the Lower Bound is defined as the maximum lower bound, $LB = max \{LB_1,LB_2\}$. 

It is worth noticing that with a valid $LB$ we can redefine a new date for the beginning of the processing time of jobs in the $WSPT-TRD$ algorithm. Particularly, in {\it Step 4},  $T_0= LB - \sum_{i \in J^2} p_{2i}$, this is valid for the jobs with $w_j^2 = 0$ or for all jobs when $\delta$ $\geq -1$.

\section{Hybrid Lagrangian Metaheuristic}
\label{Metah}

In the Hybrid Lagrangian Metaheuristic Framework proposed in this work, the Lagrangian Dual is solved by the Volume algorithm, as proposed in \cite{Barahona2000} and \cite{Nogueira}. Readers are referred to \cite{Barahona2006} and \cite{Fukuda2007} for a discussion on the Volume algorithm and its performance. In our case, the Lagrangian Relaxation incorporates heuristics as internal procedures with a focus on obtaining feasible solutions. These heuristics are based on a constructive method and a Local Search. Both procedures are executed sequentially under the Lagrangian Relaxation. 

The Volume algorithm is an extension of the subgradient algorithm, which produces a sequence of primal and dual solutions, thus being able to prove optimality. This algorithm has similar computational effort than the subgradient algorithm, and it presents similarities with the Conjugate Subgradient method \cite{Lemarechal1975}, \cite{Wolfe1975} and the Bundle method \cite{Lemarechal2001, Lemarechal1989}. A discussion of its main features and a global convergence analysis can be found in \cite{Bahiense2002}.

In the proposed algorithm the precedence constraints are dualized. The Lagrangian multipliers define a penalty to a given job allocated at a given position if the job has a large associated value, it means that it has a greater impact on the objective function, i.e., it is allocated in a disadvantageous position. This information can be used to decide when to schedule the jobs. 
Following, we describe two constructive heuristics used to obtain feasible solutions.

\subsection{Constructive Heuristics}

The heuristic H1 uses information from Lagrangian multipliers to sort jobs on M1, while H2 uses information from Lagrangian multipliers to sort jobs on M2. In the final step both heuristics use an NEH-like heuristic (\cite{Nawaz1983}) to construct the feasible solution. For a given list of jobs, the NEH function construct a solution by scheduling the jobs on the list one by one, in the best position of the partial schedule. The two heuristics are described below:\\

\textbf{H1}
\begin{description}
\item{\textit{Step 1:}} For each job on the machine 1 (\textit{$k \in J^1$}), compute $\beta_{k} = \dfrac{\sum_{j \in J^2} \lambda_{jk}}{Nprec_{k}}$.
 
Where $Nprec_{k}$ corresponds to the number of jobs in the second stage that are waiting for the completion of job $k$ on M1.

\item{\textit{Step 2:}} Get the sequence on M1 by ordering jobs in decreasing order of $\beta_{k}$.

\item{\textit{Step 3:}} Calculate the new release dates ($r_j$, $j \in J^2$) of jobs on M2.

\item{\textit{Step 4:}} Get the sequence on M2 by ordering jobs in non-decreasing order of release dates $r_j$ and apply NEH.

\end{description}

\textbf{H2}
\begin{description}
\item{\textit{Step 1:}} For each job on M2 (\textit{$j \in J^2$}), compute $\beta_{j} = \dfrac{\sum_{k \in S_j} \lambda_{jk}}{Nprec_{j}}$.
 
Where $Nprec_{j}$ corresponds to the number of precedents that job $j$ has.

\item{\textit{Step 2:}} Get the sequence on M2 by ordering jobs in increasing order of $\beta_{j}$.

\item{\textit{Step 3:}} Sequence the jobs on the M1 according the sequence on M2, respecting the precedence relations of 
cross-docking.

If there is a job on M1 without precedence relationship in the machine 2, schedule this job last.

\item{\textit{Step 4:}} Calculate the new release dates ($r_j$, $j \in J^2$) of jobs on M2, from the sequence on M1.

\item{\textit{Step 5:}} Get the sequence on M2 by ordering jobs in non-decreasing order of release dates $r_j$ and apply NEH.

\end{description}

\subsection{The Lagrangian algorithm}

Let $x$ and $y$ be solutions of subproblems X and Y, respectively, with objective function value $z$. 
Then $\lambda$ are the Lagrangian multipliers obtained by the Volume Algorithm, and $\nu$($\lambda$, $x$, $y$) is the subgradient. Let UB be the upper bound or feasible solution. The initial UB is generated by a Lagrangian metaheuristic. The algorithm is described below:

\begin{description}
\item{\textit{Step 0:}} We start with a null vector $\lambda_{jk}$ ($\lambda_{jk}$=0). Solve the Lagrangian subproblem to obtain the initial solution of 
X subproblem ($x^0$), Y subproblem ($y^0$) and objective function value ($z^0$). Let UB be the feasible solution obtained by Lagrangian metaheuristic. 
Set $\bar{x} = x^0$, $\bar{y}$=$y^0$, $\bar{z}=z^0$ and $k$=1.

\item{\textit{Step 1:}} Compute the subgradient $\nu$($\lambda_{jk-1}$, $\bar{x}$, $\bar{y}$) and $\lambda_{jk}$=$\bar{\lambda_{jk}}$ + $s$$\nu$($\lambda_{jk-1})$, 
the calculation of the step size $s$ is given by the equation (\ref{step}). Solve the Lagrangian subproblem with the new $\lambda_{jk}$ and 
let $x^k$, $y^k$ and $z^k$ be the solutions obtained. Then $\bar{x} = \alpha$$x^k$ + (1-$\alpha$)$\bar{x}$, $\bar{y} = \alpha$$y^k$ + (1-$\alpha$)$\bar{y}$ 
and $\bar{\lambda}$=$\alpha$$\lambda^k$ + (1-$\alpha$)$\bar{\lambda}$, where $\alpha$ is a number between 0 and 1, defined by convex 
combination such as defined in (\ref{cconvexa}).

\item{\textit{Step 2:}} If $z^k$ $>$ $\bar{z}$ update $\bar{\lambda}$=$\lambda_{jk}$ and $\bar{z}$=$z^k$. If $\bar{z}$ improves by 10\% since the 
last run of the Lagrangian metaheuristic then go to Step 3. Else go to Step 4.

\item{\textit{Step 3:}} Lagrangian heuristics (H1, H2):

\begin{enumerate}
\item $UB^k$ $\leftarrow$ Constructive Heuristics (H1 and H2);
\item $UB^k$ $\leftarrow$ Local Search (list);
\end{enumerate}
If $UB^k < UB$ update UB = $UB^k$ and add a cut in order to reduce the time horizon, improving the limits found (update $T_f$=UB). 

\item{\textit{Step 4:}} Stop criteria: If satisfied stop. Else let $k = k + 1$ and go to Step 1.
\end{description}

The step size $s$ and the parameter $\alpha$ used to define $\bar{x}$ and $\bar{y}$, are computed as proposed in \cite{Barahona2000} and
\cite{Fukuda2007}. First, we define $\alpha_{opt}$ as:

\begin{equation}
\alpha_{opt} = argmin \parallel \alpha\nu'^k_{(\lambda_{jk-1}, x^k, y^k)}  + (1-\alpha)\nu^k_{(\lambda_{jk-1}, \bar{x}, \bar{y})} \parallel^2 \label{cconvexa}
\end{equation}

The parameter values $\pi$, {\it MaxWaste}, {\it  factor}, $\alpha_{max}$, {\it  st}, $\alpha$, {\it yellow} and {\it green}, some of them yet to be introduced, were defined using the SPOT method, as explained in Appendix \ref{SPOT}. The parameters $\alpha_{max}$ and $\alpha$ are initially defined as 0.3034 and 0.0830 respectively, based on 
computational tests performed by SPOT. And $\alpha$ is computed as:

\begin{equation}
 \alpha={\alpha_{max} * \alpha} \hspace{0.3cm} if \hspace{0.3cm} \alpha_{opt} <0
\end{equation}

\begin{equation}
 \alpha= min\{\alpha_{opt},\alpha_{max}\} \hspace{0.3cm} if \hspace{0.3cm} \alpha_{opt} \geq0
\end{equation}

Before setting the value of the step $s$, we need to define the parameter $\pi$. Therefore,  to set the value of $\pi$ we define three types of algorithms iterations as defined in \cite{Barahona2000} and \cite{Barahona2006}.

Each iteration with no improvement is named {\it red}. Then, the parameter {\it MaxWaste} represents the maximum number of iterations without a lower bound improvement. If $z^k$ $>$ $\bar{z}$ and $\nu'^k_{(\lambda_{jk-1}, x^k, y^k)} \nu^k_{(\lambda_{jk-1}, \bar{x}, \bar{y})} < 0$, it means that a longer step in the direction to $\nu^k$ would have given a smaller value for $z^k$. Those 
iterations are denominated \textit{yellow}, otherwise, the iteration is denominated \textit{green}. At each \textit{yellow} iteration we would 
multiply $\pi$ by the {\it yellow} parameter. At each \textit{green} iteration we would multiply $\pi$ by the {\it green} parameter. After a sequence of 24 
consecutive \textit{red} iterations, we would multiply $\pi$ by {\it factor} parameter. Thus, the step size $s$ at iteration $k$ is defined as:

\begin{equation}
s = \dfrac{\pi * (st*UB - \bar{z}))}{\parallel \nu_{(\lambda_{jk-1}, \bar{x}, \bar{y})}\parallel^2} \label{step}
\end{equation}

In the “Lagrangian heuristics step” two constructive heuristics and a Local Search are executed sequentially. The Local Search is implemented based on the proposals of \cite{Arroyo2009} and \cite{Stutzle1998}. The procedure adopted is composed by \textquotedblleft swap\textquotedblright \space and \textquotedblleft insertion\textquotedblright \space moves in the sequence on machine 2 generated by H1 and H2. The former consists of interchanging all pairs of jobs. The latter consists of removing a job from its original position and inserting it on one of the $n-1$ remaining positions. The local search procedure stops when it is unable to improve the solution further.

Finally, the stop criterion is determined if one of the following criteria is satisfied:
\begin{enumerate}
 \item Maximum number of iterations, in this case 1000 iterations, or;
 \item Relative tolerance GAP defined as: $\dfrac{z^k - \bar{z}}{\bar{z}} < 0.1$, or;
 \item Null Subgradient module or negligible: $\parallel \nu_{(\lambda_{jk-1}, \bar{x}, \bar{y})}\parallel^2  \leq  0.000001$.
\end{enumerate}

%

\subsection{JB Heuristic}
\label{JB}

We compare the results of heuristics H1 and H2 with the best heuristic introduced by Chen and Lee \cite{ChenLee2009}, which is based on Johnson's algorithm. 
The heuristic JB is proposed in two steps. Firstly, for an instance of $F2|CD|C_{max}$, the authors construct an instance of $F2||C_{max}$ with $n$ jobs. The first stage is maintained unchangeable and the second stage is converted into $n$ jobs. Secondly, Johnson's algorithm is
applied to obtain a sequence in the first stage, generating a lower bound. From the sequence in the first stage, 
jobs $j \in J_{2}$ at the second stage are sequenced as soon as possible in the time horizon, respecting the completion time of its predecessors. 
We run the above algorithm for the primary problem and its reverse problem and select the best solution. As a result, an Upper Bound for the problem is generated (See \cite{ChenLee2009} for a more detailed explanation), so we calculated the relative GAP. With the objective of a fair comparison, the NEH algorithm and Local search are used to refine the results of JB in the same way than used in H1 and H2.

\subsection{CDH Heuristic}
\label{CDH}

The proposed H1 and H2 algorithms are also compared with the heuristic in Cota et al. \cite{CotaEtAll2016}, called CDH.
The CDH heuristic starts by creating a preliminary schedule in the second machine.
With this, a schedule is defined in the first machine and finally the second machine is rescheduled. The heuristic 
initially creates fictitious processing times $TF_{j} = (\sum_{i \in J^1,i \in S_j} p_{1i}) + p_{2j}$ for each job $j \in J^2$. Then, a preliminary schedule on M2 
is obtained by ordering jobs in increasing order of fictitious processing times. Then the jobs on the M1 are sequence 
according the preliminary schedule on M2, respecting the precedence relations of 
cross-docking. The sequence on M1 imposes release dates to the jobs in the second machine due to precedence relationships.
So, the ERD rule (earliest release date first) is applied to reschedule jobs in the second machine, obtaining a sequence on the M2. 
Let UB be the upper bound obtained with the heuristic and LB be the lower bound calculated as proposed in \cite{ChenSong2009}, we generate the percentage GAP for all tested instances.

\section{Generalization for parallel-docks CDC}
\label{MD}
In this section, we present a time-indexed mixed integer formulation for $F2(P)|CD|C_{max}$ based on the proposal of Cota et al. \cite{CotaEtAll2016}.
The cross-docking problem is modeled as a hybrid two-stage flow shop scheduling 
problem with identical machines and cross-docking constraints, with the objective of minimizing the makespan. 
The problem is denoted as $F2(P)|CD|C_{max}$, since $F2|CD|C_{max}$ is strongly NP-hard, as showed by Chen and Lee \cite{ChenLee2009}, it is not difficult to see that $F2(P)|CD|C_{max}$ is also strongly NP-hard.

\subsection{Time-indexed mixed integer linear programming model for $F2(P)|CD|C_{max}$}

For the parallel-docks formulation, we use the following notation:\\

\textbf{Input Parameters}
\begin{itemize}
 \item $m_1$: number of parallel processors in stage 1.
 \item $m_2$: number of parallel processors in stage 2.
  \item $n_1$: number of jobs in stage 1.
  \item $n_2$: number of jobs in stage 2.
  \item $p_{j}$: processing time of job $j$.
    \item $T$: Time horizon considered initially. In practice, we used $T = \sum_{{j} \in J_1} p_{j} + \sum_{{j} \in J_2} p_{j}$.
\end{itemize}

\textbf{Sets}
\begin{itemize}

\item $T$ = \{$0,...,T$\}, set of discrete periods considered.

\item  $J_1$ = \{$1, 2,...,n_1$\}, set of jobs in stage 1.

\item  $J_2$ = \{$1, 2,...,n_2$\}, set of jobs in stage 2.

\item $S_{j}$: a set of precedent subset jobs of $J_1$ corresponding to job $j \in J_2$, $S_{j} \subset  S = \{S_{1} ,..., S_{{n_2}} \}$.
\end{itemize}

\textbf{Decision variables}
\begin{itemize}
 \item $C_{max}$: makespan.
  \item $x_{jt}=1$, if job $j$ starts to be processed at time $t$, 0, otherwise.
\end{itemize}
The multiple dock formulation is defined as:

\newpage
\begin{fleqn}
\begin{equation}
\hspace{2 cm} {Minimize}   \hspace{0.2cm} {C_{max} \label{foeqmd}} 
\end{equation}
\end{fleqn}

subject to\\
\begin{eqnarray}
\sum_{t=0}^{T-p_j} x_{jt} = 1, & & \forall j \in J_1, \label{eqmd01} \\
\sum_{t=0}^{T-p_j} x_{jt} = 1, & & \forall j \in J_2, \label{eqmd02} \\
\sum_{j \in J_1} \sum_{s=max{(0 ; t-p_j+1)}}^t x_{js} \leq m_1, & & \forall t \in T, \label{eqmd03}\\
\sum_{j \in J_2} \sum_{s=max{(0 ; t-p_j+1)}}^t x_{js} \leq m_2, & & \forall t \in T, \label{eqmd04}\\
\sum_{t=0}^{T-p_j}tx_{jt} \geq \sum_{t=0}^{T}(t+p_i)x_{it}, & & \forall j \in J_2, \forall i \in S_j, \label{eqmd05}\\
C_{max} \geq p_j + \sum_{t=0}^{T-p_j}tx_{jt}, & & \forall j \in J_2, \label{eqmd06}\\
x_{jt} \in \{0,1\}, & &  \forall j \in J_1 \wedge \forall j \in J_2,  \forall t \in T, \label{eqmd07}\\
C_{max} \geq 0. \label{eqmd08}
\end{eqnarray}

The objective function remains as the minimization of the makespan. Constraints (\ref{eqmd01}) and (\ref{eqmd02})
ensure that in every stage, each job should start its processing on one, and only one, time slot within the time horizon $T$. 
Constraints set (\ref{eqmd03}) and (\ref{eqmd04}) indicates that 
the number of jobs that can be process simultaneously is less than or equal to the number of existing parallel processors in the stage considered. 
The constraint set (\ref{eqmd05}) represents the cross-docking constraints, in which the release date of each job $j \in J^2$, must be greater than or equal to the completion date of each 
task belonging to its precedence set. Besides that, the constraint set ensures that each job $j \in J^2$, starts processing only after its release date. 
Constraints (\ref{eqmd06}) compute the makespan, analyzing the maximum completion time. Constraints (\ref{eqmd07}) and (\ref{eqmd08}) specify the domains of each decision variable.

\subsection{Lagrangian Relaxation for $F2(P)|CD|C_{max}$}

Let  $L_{MD}(\lambda$) be the relaxation of multiple dock formulation, when the cross-docking constraints (\ref{eqmd05}) are dualized and its violation 
penalized in the objective function. We define the subproblem $L_{MD}$($\lambda$) as:

\begin{fleqn}
\begin{equation}
 L_{MD}(\lambda) = Minimize \hspace{0.2cm} {C_{max}} - \sum_{j \in J_2} \sum_{i \in S_j} \lambda_{ij} (\sum_{t=0}^{T-p_j}tx_{jt} - \sum_{t=0}^{T}(t+p_i)x_{it})) \label{eqlgmd}
\end{equation}
\end{fleqn}
s.t. (\ref{eqmd01}) -- (\ref{eqmd04}), (\ref{eqmd06}) -- (\ref{eqmd08}), $\lambda_{ij} \geq 0$.

Just as we did for the 2-dock case, we are able to uncouple the problem into two new scheduling problems, one for each stage. 

\textbf{Subproblem 1}

\begin{equation}
 L'_{MD}(\lambda) = Min \hspace{0.2cm} \sum_{j \in J_2} \sum_{i \in S_j} \lambda_{ij} \sum_{t=0}^{T}(t+p_i)x_{it}  \label{eqsublg1}
\end{equation}
subject to\\
\begin{eqnarray}
\sum_{t=0}^{T-p_j} x_{jt} = 1, & & \forall j \in J_1, \label{eqlg101} \\
\sum_{j \in J_1} \sum_{s=max{(0 ; t-p_j+1)}}^t x_{js} \leq m_1, & & \forall t \in T, \label{eqlg102}\\
x_{jt} \in \{0,1\}, & &  \forall j \in J_1, \forall t \in T, \label{eqlg103}\\
\lambda_{ij}  \geq 0, & & \forall j \in J_2, i \in S_j. \label{eqlg104}
\end{eqnarray}

\textbf{Subproblem 2}

\begin{equation}
 L''_{MD}(\lambda) = Min \hspace{0.2cm} {C_{max}} - \sum_{j \in J_2} \sum_{i \in S_j} \lambda_{ij} \sum_{t=0}^{T-p_j}tx_{jt}  \label{eqsublg2}
\end{equation}
subject to\\
\begin{eqnarray}
\sum_{t=0}^{T-p_j} x_{jt} = 1, & & \forall j \in J_2, \label{eqlg201} \\
\sum_{j \in J_2} \sum_{s=max{(0 ; t-p_j+1)}}^t x_{js} \leq m_2, & & \forall t \in T, \label{eqlg202}\\
C_{max} \geq p_j + \sum_{t=0}^{T-p_j}tx_{jt}, & & \forall j \in J_2, \label{eqlg203}\\
x_{jt} \in \{0,1\}, & &  \forall j \in J_2, \forall t \in T, \label{eqlg204}\\
\lambda_{ij}  \geq 0, & & \forall j \in J_2, i \in S_j, \label{eqlg206}\\
C_{max} \geq 0. \label{eqlg207}
\end{eqnarray}

\subsection{Constructive Heuristic for parallel-dock CDC}

The first subproblem can be defined as $Pm||\sum C_j W_j$, while the second as $Pm||C_{max} - \sum I_j W_j$, therefore, it is not difficult to see that both subproblem are strongly NP-hard (see \cite{Pinedo2008}). In this case, we focus on improving the problem's upper bounds, thus, in each iteration, we solve the linear relaxations of the subproblems, obtain integer solutions and then compute their Lagrangian multipliers. The logic of the procedure is described in Algorithm \ref{algMD}, maintaining the same reasoning proposed in Section 4. 

For comparison purposes, we limited the computational time in 60 seconds, and we compute the LB precisely as proposed by Chen and Song \cite{ChenSong2009}.

\begin{algorithm}
\caption{Algorithm MD}\label{algMD}

\begin{algorithmic}[1]
\State Set Lagrangian multipliers as null.
\State Compute a Lower Bound (LB), as proposed by Chen and Song \cite{ChenSong2009}. 

\State Compute an initial Upper Bound (UB).

\begin{align*}
UB = &\sum_{i \in J_1} p_i/m_1 - max_{i \in J_1}\hspace{0.1cm} p_i/m_1 + max_{i \in J_1}\hspace{0.1cm} p_i + \\
&+\sum_{j \in J_2} p_j/m_2 - max_{j \in J_2}\hspace{0.1cm} p_j/m_2 + max_{j \in J_2}\hspace{0.1cm} p_j.
\end{align*}

\While{Time < TimeLimit}
  \State Solve the Linear Relaxation.
  \State \ \ \ \ Allocate the job as long as its variable is closer to 1.
  \State \ \ \ \ Calculate the Lagrangian multipliers.
  \State \ \ \ \ \ Compute the weights $\beta_{1}$ and $\beta_{2}$ for the heuristics H1 and H2 (the ordered weights define the entry list for H1 and H2).
  \State \ \ \ \ Apply Heuristics H1 and H2.
  \State \ \ \ \ Apply a First-Best Local Search.
  \State \ \ \ \ \textbf{if}  UB has improved
  \State \ \ \ \ \ \ Update the UB and add a cut to reduce the time horizon.
  \State \ \ \ \ \textbf{else} 
    \State \ \ \ \ \ \  \textbf{end algorithm}.
  \State \ \ \ \ \ \  \textbf{end while}.
  \State \ \ \ \ \textbf{end if}
\State \textbf{end while}
\end{algorithmic}
\end{algorithm}

\newpage
\section{Computational Experiments}
\label{CE}

\subsection{Results for the 2-dock problem}

To investigate the performance of the Complete Model and Hybrid Lagrangian Metaheuristic Framework, artificial instances are generated varying the processing time of jobs and the number of jobs on machines 1 and 2, according to Chen and Lee \cite{ChenLee2009}. 
We ran the tests on a single thread in the Intel Xeon-Silver 4110 (2.1GHz/8-core/85W) with 64 GB memory and Linux operational system. The programming language used is C++ with the optimization software CPLEX 12.4.

\subsubsection{Instances Generation}

The instances used in this work are generated through the software MATLAB, following the description found in Chen and Lee~\citep{ChenLee2009}. 
Table \ref{tabelainst} presents a summary of the generated instances and its characteristics. 
It is divided into two groups, each row informs the features of a sub-group of instances. $n$ and $m$ indicate the number of jobs in each stage. $(NP)$ the discrete uniform distribution to select the number predecessors of jobs $j \in J^2$ and $(TP)$ the distribution to generate all processing times.

In each group of instances, $n$ is fixed and the number $m$ is randomly chosen from the range of $[0.6n,0.8n, 1.0n, 1.2n, 1.4n]$. For each $n$, five instances with different $m$ values are considered. And for each pair  ($n$, $m$), we generate 10 different problems, therefore the benchmark consists of 500 instances, 50 instances per row of Table \ref{tabelainst} \footnote{all instances and codes are available at \url{https://github.com/GabrielaBragaFonseca/Cross-docking-Problems}}.

\begin{table}
 \centering
\caption{Summary of the artificial benchmark.} 
  \begin{tabular}{c|cccc}
   \hline
   Group &  Jobs &  Jobs & NP    & TP     \\
          & M1(n) &  M2(m) &       &      \\ \hline \hline
          & 5      & 3-4-5-6-7 & U(1,4)     & U(1,10) \\
          & 10      & 6-8-10-12-14 & U(1,9)     & U(1,10)  \\
    1	  & 20    & 12-16-20-24-28 & U(1,19)     & U(1,10)  \\
    	  & 40    & 24-32-40-48-56 & U(1,39)     & U(1,10) \\
	  & 60    & 36-48-60-72-84 & U(1,59)     & U(1,10) \\\hline\hline	  
          & 5      & 3-4-5-6-7 & U(1,4)     & U(10,100) \\
          & 10      & 6-8-10-12-14 & U(1,9)     & U(10,100)  \\
    2	  & 20    & 12-16-20-24-28 & U(1,19)     & U(10,100)  \\
    	  & 40    & 24-32-40-48-56 & U(1,39)     & U(10,100)  \\
	  & 60    & 36-48-60-72-84 & U(1,59)     & U(10,100) \\\hline\hline	  
 
  \end{tabular} 
  \label{tabelainst}
\end{table} 

\subsubsection{Computational Results}

In Table \ref{tabelaresultados2} we report computational results for the Complete Model (MIP) and heuristics JB, CDH, H1, and H2. Columns $n$ and $m$ show respectively the number of jobs in the first and second machine. For each instance size ($n$, $m$) we report the best result and the average of 10 cases tested. We also presented the subgroup's  average results.

The GAP(\%) is computed as $\dfrac{(\mbox{Upper Bound - Lower Bound)}}{\mbox{Upper Bound}}$ and the column 
T(s) refers to CPU time expended to solve the problem in seconds. The run time is limited to one hour of CPU time 
(3.600 seconds), and results are depicted in Table \ref{tabelaresultados2} are registered. The dash (-) means that 
the corresponding value is not found. Finally, after comparing and analyzing the MIP, JB, CDH, H1 and H2 the best 
results obtained are highlighted in the Table.

The GAP is calculated considering the best lower bound found when comparing the lower bounds of Complete Model and 
Hybrid Lagrangian Metaheuristic Framework. It is worth highlighting that the constructive heuristic JB is proposed 
by Chen and Lee \cite{ChenLee2009} and CDH heuristic is proposed by Cota et al. \cite{CotaEtAll2016}.

The MIP results show that the proposed model is efficient to solve to optimality only small instances. 
Furthermore, for the larger instances the model found no solution. These results expose the difficulty of solving time-indexed models. 
In the time-indexed problems, the number of variables is proportional to the time horizon, and the higher the number of jobs the greater the horizon of time to sequence them,  increasing the computational effort to solve these problems. This fact justifies the proposal of the hybrid constructive heuristics to solve instances with a greater number of  jobs. 

Results show that H1 and H2 performs better than heuristics JB and CDH. Based on the Table 
\ref{tabelaresultados2}, the heuristic H2 showed GAP results more efficients. Comparing the average GAPs we can 
observe that the heuristic H2 has its better performance when $n < m$ for small and large instances. When $n > m$, 
the heuristic H2 gets better for small cases, while H1 presents better results for large instances. 
The same behavior is noticed when $n = m$, H2 is best for small cases, and H1 for larger ones. The heuristic JB 
had no best GAP in any instance. As for the CDH heuristic, the GAPs found were high, which was already expected, since 
the heuristic was developed for a more general machine's scenario.

\begin{table}[htbp]
  \caption{Computational results for complete model and constructive heuristics.} 
 \begin{center}
  \scalebox{0.5}{
  \begin{tabular}{ccccc|cc|cc|cc|cc||cc|cc|cc|cc|cc}
   \hline
   &       & \multicolumn{10}{c}{Group 1 - Processing time  [1, 10]} & \multicolumn{10}{c}{Group 2 - Processing time  [10, 100]} \\
   \hline
   \multicolumn{1}{c}{n} & \multicolumn{1}{c}{m} &       & \multicolumn{2}{c}{MIP} &  \multicolumn{2}{c}{JB}  &  \multicolumn{2}{c}{CDH} &  \multicolumn{2}{c}{H1} &  \multicolumn{2}{c||}{H2} &  \multicolumn{2}{c}{MIP} &  \multicolumn{2}{c}{JB}  &  \multicolumn{2}{c}{CDH} &  \multicolumn{2}{c}{H1} &  \multicolumn{2}{c}{H2}  \\
			 &       &       & GAP & T(s) 		& GAP &  T(s)              &  GAP   & T(s) 		&  GAP   & T(s)               &  GAP   & T(s) & GAP & T(s) 		& GAP &  T(s)              &  GAP   & T(s) 		&  GAP   & T(s) &  GAP   & T(s)\\
    \hline
    \multicolumn{1}{c}{5} & \multicolumn{1}{c}{3} & \multicolumn{1}{c}{Best} &  0.0\%  & 0.1 		& \textbf{0.0\%}  & 0.0		   & 15.2\%  & 0.0 & \textbf{0.0\%}  & 0.0		& \textbf{0.0\%}  & 0.0     						& 0.0\%  & 6.7 		& \textbf{0.0\%}  & 0.0	& 7.8\%  & 0.0    & \textbf{0.0\%}  & 0.0 		& \textbf{0.0\%}  & 0.0  \\
    &       & \multicolumn{1}{c}{Average} 				     & 0.0\%  & 0.2 		& 6.8\%  & 0.0  	   & 26.3\%  & 0.0 & \textbf{4.8\%}  & 0.0 		& \textbf{4.8\%}  & 0.0   	 									& 0.0\%  & 21.0 		& 2.4\%  & 0.0		   & 28.8\%  & 0.0 & \textbf{2.3\%}  & 0.0 		& \textbf{2.3\%}  & 0.0  \\  
           & \multicolumn{1}{c}{4} & \multicolumn{1}{c}{Best} &  0.0\%  & 0.1 		& \textbf{0.0\%}  & 0.0		  & 2.7\%  & 0.0	 & \textbf{0.0\%}  & 0.0		& \textbf{0.0\%}  & 0.0     									& 0.0\%  & 3.9 		& \textbf{0.0\%}  & 0.0		   & 3.0\%  & 0.0 & \textbf{0.0\%}  & 0.0 		& \textbf{0.0\%}  & 0.0 \\ 
   &       & \multicolumn{1}{c}{Average} 				     &  0.0\%  & 0.3 		& 7.7\%  & 0.0  	   & 22.3\%  & 0.0	& \textbf{5.9\%}  & 0.0 		& \textbf{5.9\%}  & 0.0     									& 0.0\%  & 354.9 		& 7.5\%  & 0.0		   & 24.2\%  & 0.0 & \textbf{6.2\%}  & 0.0		& \textbf{6.2\%}  & 0.0  \\  
           & \multicolumn{1}{c}{5} & \multicolumn{1}{c}{Best} &  0.0\%  & 0.3 & 4.5\%  & 0.0		   & 17.1\%  & 0.0 & \textbf{0.0\%}  & 0.0		& \textbf{0.0\%}  & 0.0     													& 0.0\%  & 32.6 		& 3.8\%  & 0.0		   & 19.3\%  & 0.0 & \textbf{0.0\%}  & 0.0 		& \textbf{0.0\%}  & 0.0  \\
    &       & \multicolumn{1}{c}{Average} 				     & 0.0\%  & 0.5 		& 9.8\%  & 0.0  	  & 32.8\%  & 0.0	 & \textbf{2.5\%}  & 0.0 		& 3.3\%  & 0.0     										& 0.8\%  & 1570.1 		& 9.9\%  & 0.0		   & 33.4\%  & 0.0 & \textbf{2.6\%}  & 0.0		& \textbf{2.6\%}  & 0.0  \\  
           & \multicolumn{1}{c}{6} & \multicolumn{1}{c}{Best} &  0.0\%  & 0.1  & 4.8\%  & 0.0	   & 13.0\%  & 0.0	& \textbf{0.0\%}  & 0.0 & \textbf{0.0\%}  & 0.0     														& 0.0\%  & 80.9     & 5.7\%  & 0.0		  & 13.5\%  & 0.0 & \textbf{0.0\%}  & 0.0 		& \textbf{0.0\%}  & 0.0  \\
    &       & \multicolumn{1}{c}{Average} 		      &  0.0\%  & 1.3  & 16.6\%  & 0.0  	   & 29.4\%  & 0.0	& \textbf{5.2\%}  & 0.0 & 5.5\%  & 0.0     																& 11.8\%  & 2388.3  & 16.1\%  & 0.0		   & 30.7\%  & 0.0 & \textbf{4.8\%}  & 0.0		& 6.2\%  & 0.0  \\  
           & \multicolumn{1}{c}{7} & \multicolumn{1}{c}{Best} &  0.0\%  & 0.2  & \textbf{0.0\%}  & 0.0	   & 16.7\%  & 0.0	& \textbf{0.0\%}  & 0.0 & \textbf{0.0\%}  & 0.0     													& 0.0\%  & 94.9          & \textbf{0.0\%}  & 0.0		  & 16.6\%  & 0.0 & \textbf{0.0\%}  & 0.0 		& \textbf{0.0\%}  & 0.0   \\
    &       & \multicolumn{1}{c}{Average} 		      &  0.0\%  & 3.7  & 14.1\%  & 0.0  	   & 36.4\%  & 0.0	& 6.3\%  & 0.0 & \textbf{5.4\%}  & 0.0     	   															& 13.1\%  & 2882.5       & 13.9\%  & 0.0		   & 38.0\%  & 0.0 & \textbf{5.1\%}  & 0.0		& 6.7\%  & 0.0  \\ \hdashline
    \multicolumn{3}{c}{Subgroup Average} & 0.0\% & 1.2 & 11.0\% & 0.0 & 29.4\%  & 0.0 & \textbf{4.9\%} & 0.0 & 5.0\% & 0.0 																						& 5.1\% & 1443.4  & 10.0\% & 0.0 & 31.0\%  & 0.0 & \textbf{4.2\%} & 0.0 & 4.8\% & 0.0  \\\hline
    \multicolumn{1}{c}{10} & \multicolumn{1}{c}{6} & \multicolumn{1}{c}{Best} &  0.0\%  & 5.1 		& \textbf{0.0\%}  & 0.0		   & 7.5\%  & 0.0 & \textbf{0.0\%}  & 0.0		& \textbf{0.0\%}  & 0.0    									& 28.7\%  & 2801.1 		& 0.7\%  & 0.0	   & 8.9\%  & 0.0 & \textbf{0.0\%}  & 0.0 		& \textbf{0.0\%}  & 0.0  \\
    &       & \multicolumn{1}{c}{Average} 				     &  0.0\%  & 20.2 		& 11.8\%  & 0.0  	   & 24.8\%  & 0.0	& 1.8\%  & 0.0 		& \textbf{1.7\%}  & 0.0     	  											& 28.2\%  & 3509.9 	        & 11.3\%  & 0.0		  & 23.0\%  & 0.0 & \textbf{1.4\%}  & 0.0 		& \textbf{1.4\%}  & 0.0  \\  
          & \multicolumn{1}{c}{8} & \multicolumn{1}{c}{Best} & 0.0\%  & 10.2 		& 2.7\%  & 0.0		  & 9.6\%  & 0.0 & \textbf{0.0\%}  & 0.0		& \textbf{0.0\%}  & 0.0     														& 23.6\%  & 3597.6 		& 2.3\%  & 0.0	  & 15.0\%  & 0.0  & \textbf{0.0\%}  & 0.0 		& \textbf{0.0\%}  & 0.0 \\
    &	    & \multicolumn{1}{c}{Average} 		       &  0.0\%  & 167.5 		& 13.3\%  & 0.0  	   & 25.8\%  & 0.0	& 6.0\%  & 0.0 		& \textbf{5.4\%}  & 0.0     													& 36.7\%  & 3600.0 	        & 17.0\%  & 0.0		  & 30.5\%  & 0.0 & \textbf{3.7\%}  & 0.0 		& 4.2\%  & 0.0  \\  
          & \multicolumn{1}{c}{10} & \multicolumn{1}{c}{Best} & 0.0\%  & 14.1 		& 1.4\%  & 0.0		   & 14.8\%  & 0.0	& \textbf{0.0\%}  & 0.0		& 2.2\%  & 0.0     														& 35.7\%  & 3600.0 		& 2.2\%  & 0.0	   & 26.4\%  & 0.0 & 4.6\%  & 0.0 		& \textbf{1.0\%}  & 0.0  \\
    &	    & \multicolumn{1}{c}{Average} 		       &  0.0\%  & 369.6 		& 15.9\%  & 0.0  	   & 28.8\%  & 0.0	& 6.5\%  & 0.0 		& \textbf{6.4\%}  & 0.0     														& 45.2\%  & 3600.0 	        & 18.4\%  & 0.0		   & 36.8\%  & 0.0 & \textbf{7.6\%}  & 0.0 		& 7.9\%  & 0.0  \\  
          & \multicolumn{1}{c}{12} & \multicolumn{1}{c}{Best}  & 0.0\%  & 42.1 		& 3.9\%  & 0.0		   & 17.6\%  & 0.0	& \textbf{0.0\%}  & 0.0		& \textbf{0.0\%}  & 0.0     													& 51.7\%  & 3600.0 		& 3.1\%  & 0.0	   & 19.3\%  & 0.0 & \textbf{0.0}\%  & 0.0 		& \textbf{0.0}\%  & 0.0  \\
    &	    & \multicolumn{1}{c}{Average} 		       &  5.0\%  & 2440.5		& 14.0\%  & 0.0  	  & 41.0\%  & 0.0	 & 5.4\%  & 0.0 		& \textbf{5.3\%}  & 0.0     														& 65.2\%  & 3600.0 	        & 13.9\%  & 0.0		   & 40.8\%  & 0.0 & \textbf{5.4}\%  & 0.0 		& \textbf{5.4\%}  & 0.0  \\  
          & \multicolumn{1}{c}{14} & \multicolumn{1}{c}{Best}  & 0.0\%  & 21.9 		& 1.3\%  & 0.0		 & 27.9\%  & 0.0  & \textbf{0.0\%}  & 0.0		& 1.0\%  & 0.0     																& 53.2\%  & 3600.0 		& 2.0\%  & 0.0	   & 28.9\%  & 0.0 & \textbf{0.0\%}  & 0.0 		& \textbf{0.0\%}  & 0.0  \\
    &	    & \multicolumn{1}{c}{Average} 		       &  4.6\%  & 3141.6		& 10.8\%  & 0.0  	  & 45.5\%  & 0.0 & \textbf{7.0\%}  & 0.0 		& 8.1\%  & 0.0     															& 73.6\%  & 3600.0 	        & 11.4\%  & 0.0		  & 45.2\%  & 0.0 & 6.7\%  & 0.0 		& \textbf{5.8\%}  & 0.0  \\\hdashline  
    \multicolumn{3}{c}{Subgroup Average} & 1.9\% & 1227.9 & 13.1\% & 0.0 & 33.2\%  & 0.0	& \textbf{5.3\%} & 0.0 & 5.4\% & 0.0 																						& 49.8\% & 3582.0 & 14.4\% & 0.0 & 35.2\%  & 0.0 & 5.0\% & 0.0 & \textbf{4.9\%} & 0.0 \\\hline
    \multicolumn{1}{c}{20} & \multicolumn{1}{c}{12} & \multicolumn{1}{c}{Best} &  1.5\%  & 3600.0 		& 10.0\%  & 0.0		    & 12.4\%  & 0.0 & 1.5\%  & 0.0		& \textbf{0.7\%}  & 0.0          											& -  & 3600.0 		& 10.2\%  & 0.0	   & 13.7\%  & 0.0 & \textbf{0.6\%}  & 0.0 		& \textbf{0.6\%}  & 0.0  \\
    &       & \multicolumn{1}{c}{Average} 				     &  14.8\%  & 3600.0  		& 18.1\%  & 0.0  	  & 26.1\%  & 0.0 & \textbf{6.6\%}  & 0.0 		& 7.1\%  & 0.0     													& -  & 3600.0 	        & 17.9\%  & 0.0		    & 26.0\%  & 0.0 & \textbf{6.4\%}  & 0.0 		& 6.5\%  & 0.0 \\  
        & \multicolumn{1}{c}{16} & \multicolumn{1}{c}{Best}  & 14.5\%  & 3600.0  		& 16.4\%  & 0.0		   & 21.4\%  & 0.0	& \textbf{1.6\%}  & 0.0		& 5.4\%  & 0.0     														& -  & 3600.0 		& 17.4\%  & 0.0	   & 22.7\%  & 0.0 & \textbf{2.3\%}  & 0.0 		& 4.2\%  & 0.0  \\
    &	    & \multicolumn{1}{c}{Average} 		       &  32.4\%  & 3600.0 		& 23.0\%  & 0.0  	   & 32.6\%  & 0.0	& \textbf{9.8\%}  & 0.0 		& 10.5\%  & 0.0     														& -  & 3600.0 	        & 23.3\%  & 0.0		   & 32.6\%  & 0.0 & \textbf{8.2\%}  & 0.0 		& 10.4\%  & 0.0  \\  
          & \multicolumn{1}{c}{20} & \multicolumn{1}{c}{Best}  & 24.5\%  & 3600.0  		& 14.3\%  & 0.0		   & 31.9\%  & 0.0 & \textbf{2.5\%}  & 0.0		& 4.4\%  & 0.0     															& -  & 3600.0 		& 15.5\%  & 0.0	   & 31.2\%  & 0.0 & 3.7\%  & 0.0 		& \textbf{1.3\%}  & 0.0  \\
    &	    & \multicolumn{1}{c}{Average} 		       &  36.0\%  & 3600.0 		& 24.2\%  & 0.0  	   & 38.5\%  & 0.0	& 11.1\%  & 0.0 		& \textbf{10.9\%}  & 0.0     														& -  & 3600.0 	        & 24.8\%  & 0.0		  & 38.0\%  & 0.0 & 11.4\%  & 0.0 		& \textbf{10.3\%}  & 0.0  \\  
          & \multicolumn{1}{c}{24} & \multicolumn{1}{c}{Best}  & 19.9\%  & 3600.0  		& 2.4\%  & 0.0		  & 35.0\%  & 0.0	 & \textbf{0.0\%}  & 0.0		& 2.6\%  & 0.0     														& -  & 3600.0 		& 1.5\%  & 0.0	   & 35.6\%  & 0.0 & \textbf{0.6\%}  & 0.0 		& 5.2\%  & 0.0  \\
    &	    & \multicolumn{1}{c}{Average} 		       &  38.4\%  & 3600.0 		& 21.2\%  & 0.0  	   & 44.4\%  & 0.0	& \textbf{12.8\%}  & 0.0 		& 14.0\%  & 0.0     														& -  & 3600.0 	        & 21.7\%  & 0.0		   & 44.1\%  & 0.0 & \textbf{12.2\%}  & 0.0 		& 13.0\%  & 0.0 \\  
           & \multicolumn{1}{c}{28} & \multicolumn{1}{c}{Best}  & 19.4\%  & 3600.0  		& 1.6\%  & 0.0		   & 39.2\%  & 0.0 & \textbf{1.1\%}  & 0.0		& 4.9\%  & 0.0     															& -  & 3600.0 		& 3.1\%  & 0.0	   & 41.0\%  & 0.0 & \textbf{0.9\%}  & 0.0 		& \textbf{0.9\%}  & 0.0  \\
    &	    & \multicolumn{1}{c}{Average} 		       &  36.5\%  & 3600.0 		& 16.7\%  & 0.0  	  & 49.3\%  & 0.0 & 15.2\%  & 0.0 		& \textbf{15.0\%}  & 0.0     															& -  & 3600.0 	        & 17.4\%  & 0.0		   & 48.9\%  & 0.0 & 14.4\%  & 0.0 		& \textbf{13.5\%}  & 0.0 \\\hdashline
    \multicolumn{3}{c}{Subgroup Average} & 31.6\% & 3600.0 & 20.6\% & 0.0 & 38.2\%  & 0.0	& \textbf{11.1\%} & 0.0 & 11.5\% & 0.0 																								& $\ast$ & 3600.0 & 21.0\% & 0.0 & 37.9\%  & 0.0 & \textbf{10.5\%} & 0.0 & 10.7\% & 0.0  \\\hline
    \multicolumn{1}{c}{40} & \multicolumn{1}{c}{24} & \multicolumn{1}{c}{Best} &  40.2\%  & 3600.0 		& 15.7\%  & 0.0		   & 20.9\%  & 0.0 & \textbf{3.9\%}  & 0.0		& 4.5\%  & 0.0          													& -  & 3600.0 		& 2.3\%  & 0.0	   & 50.0\%  & 0.0 & \textbf{0.0\%}  & 0.0 		& 1.8\%  & 0.0 \\
    &       & \multicolumn{1}{c}{Average} 				     &  46.3\%  & 3600.0  		& 20.7\%  & 0.0  	  & 26.2\%  & 0.0 & 9.6\%  & 0.0 		& \textbf{9.0\%}  & 0.0     														& -  & 3600.0 	        & 8.6\%  & 0.0		   & 56.0\%  & 0.0 & 6.2\%  & 0.0 		& \textbf{5.6\%}  & 0.0  \\  
           & \multicolumn{1}{c}{32} & \multicolumn{1}{c}{Best}  & 52.6\%  & 3600.0  		& 21.3\%  & 0.0		   & 25.2\%  & 0.0	& 9.0\%  & 0.0		& \textbf{8.3\%}  & 0.0     															& -  & 3600.0 		& 22.4\%  & 0.0	   & 27.4\%  & 0.0 & \textbf{6.0\%}  & 0.0 		& \textbf{6.0\%}  & 0.0  \\
    &	    & \multicolumn{1}{c}{Average} 		       &  57.0\%  & 3600.0 		& 26.5\%  & 0.0  	  & 31.8\%  & 0.0	 & 12.8\%  & 0.0 		& \textbf{12.5\%}  & 0.0     																& -  & 3600.0 	        & 28.5\%  & 0.0		  & 32.7\%  & 0.0  & 11.8\%  & 0.0 		& \textbf{11.0\%}  & 0.0 \\  
           & \multicolumn{1}{c}{40} & \multicolumn{1}{c}{Best}  & 57.3\%  & 3600.0  		& 26.1\%  & 0.0		  & 31.7\%  & 0.0	 & 9.1\%  & 0.0		& \textbf{8.3\%}  & 0.0     															& -  & 3600.0 		& 24.4\%  & 0.0	   & 32.6\%  & 0.0 & \textbf{9.0\%}  & 0.0 		& 10.4\%  & 0.0 \\
    &	    & \multicolumn{1}{c}{Average} 		       &  69.2\%  & 3600.0 		& 31.1\%  & 0.0  	   & 36.1\%  & 0.0	& 12.5\%  & 0.0 		& \textbf{12.3\%}  & 0.0     															& -  & 3600.0 	        & 32.4\%  & 0.0		    & 37.7\%  & 0.0 & \textbf{13.3\%}  & 0.0 		& 13.6\%  & 0.0\\  
            & \multicolumn{1}{c}{48} & \multicolumn{1}{c}{Best}  & -  & 3600.0  		& 23.9\%  & 0.0		  & 41.6\%  & 0.0	 & \textbf{9.5\%}  & 0.0		& 11.6\%  & 0.0     															& -  & 3600.0 		& 23.0\%  & 0.0	  & 41.7\%  & 0.0 & \textbf{10.8\%}  & 0.0 		& 11.6\%  & 0.0  \\
    &	    & \multicolumn{1}{c}{Average} 		       & -  & 3600.0 		        & 29.6\%  & 0.0  	   & 45.2\%  & 0.0 & \textbf{16.6\%}  & 0.0 		& 16.8\%  & 0.0     																& -  & 3600.0 	        & 29.1\%  & 0.0		   & 45.3\%  & 0.0 & 16.6\%  & 0.0 		& \textbf{16.1\%}  & 0.0  \\  
            & \multicolumn{1}{c}{56} & \multicolumn{1}{c}{Best}  & -  & 3600.0  		& 16.8\%  & 0.0		  & 44.4\%  & 0.0	 & \textbf{13.2\%}  & 0.0		& 15.5\%  & 0.0     															& -  & 3600.0 		& 18.4\%  & 0.0	   & 44.8\%  & 0.0 & \textbf{13.8\%}  & 0.0 		& 14.8\%  & 0.0 \\
    &	    & \multicolumn{1}{c}{Average} 		       & -  & 3600.0 		        & 23.3\%  & 0.0  	  & 48.3\%  & 0.0	 & 19.7\%  & 0.0 		& \textbf{19.3\%}  & 0.0     																& -  & 3600.0 	        & 23.3\%  & 0.0		   & 48.3\%  & 0.0 & 18.7\%  & 0.0 		& \textbf{18.1\%}  & 0.0  \\\hdashline  
   \multicolumn{3}{c}{Subgroup Average} & $\ast$ & 3600.0 & 26.2\% & 0.0 & 37.5\%  & 0.0	& 14.2\% & 0.0 & \textbf{14.0\%} & 0.0 																								& $\ast$ & 3600.0 & 24.4\% & 0.0 & 44.0\%  & 0.0 & 13.3\% & 0.0 & \textbf{12.9\%} & 0.0  \\\hline
   \multicolumn{1}{c}{60} & \multicolumn{1}{c}{36} & \multicolumn{1}{c}{Best} &  -  & 3600.0 		& 16.9\%  & 0.0		  & 20.9\%  & 0.0 & 5.5\%  & 0.0		& \textbf{4.5\%}  & 0.0           															& -  & 3600.0 		& 17.8\%  & 0.0	   & 21.4\%  & 0.0 & 5.9\%  & 0.0 		& \textbf{3.4\%}  & 0.0 \\
    &       & \multicolumn{1}{c}{Average} 				       &  -  & 3600.0  		& 22.8\%  & 0.0  	  & 26.7\%  & 0.0	 & 8.6\%  & 0.0 		& \textbf{8.4\%}  & 0.0     																& -  & 3600.0 	        & 22.7\%  & 0.0		   & 26.6\%  & 0.0 & 8.4\%  & 0.0 		& \textbf{7.6\%}  & 0.0  \\  
            & \multicolumn{1}{c}{48} & \multicolumn{1}{c}{Best}  & -  & 3600.0  		& 26.0\%  & 0.0		   & 30.4\%  & 0.0	& \textbf{9.6\%}  & 0.0		& 10.7\%  & 0.0     																	& -  & 3600.0 		& 26.5\%  & 0.0	  & 30.6\%  & 0.0 & 9.4\%  & 0.0 		& \textbf{8.5\%}  & 0.0  \\
    &	    & \multicolumn{1}{c}{Average} 		       & -  & 3600.0 		        & 29.9\%  & 0.0  	   & 33.0\%  & 0.0	& \textbf{12.0\%}  & 0.0 		& 12.8\%  & 0.0     																	& -  & 3600.0 	        & 29.6\%  & 0.0		  & 32.9\%  & 0.0 & \textbf{11.2\%}  & 0.0 		& 12.0\%  & 0.0   \\  
            & \multicolumn{1}{c}{60} & \multicolumn{1}{c}{Best}  & -  & 3600.0  		& 26.6\%  & 0.0		   & 35.3\%  & 0.0	& 14.2\%  & 0.0		& \textbf{13.9\%}  & 0.0     																	& -  & 3600.0 		& 28.0\%  & 0.0	   & 36.1\%  & 0.0 & 13.3\%  & 0.0 		& \textbf{12.9\%}  & 0.0  \\
    &	    & \multicolumn{1}{c}{Average} 		       & -  & 3600.0 		        & 33.5\%  & 0.0  	   & 40.6\%  & 0.0	& \textbf{18.2\%}  & 0.0 		& 18.3\%  & 0.0     																	& -  & 3600.0 	        & 33.6\%  & 0.0		   & 40.7\%  & 0.0 & \textbf{16.9\%}  & 0.0 		& 17.0\%  & 0.0  \\  
             & \multicolumn{1}{c}{72} & \multicolumn{1}{c}{Best}  & -  & 3600.0  		& 25.7\%  & 0.0		   & 42.6\%  & 0.0 & 15.5\%  & 0.0		& \textbf{14.8\%}  & 0.0     																	& -  & 3600.0 		& 26.0\%  & 0.0	  & 42.8\%  & 0.0 & \textbf{13.5\%}  & 0.0 		& \textbf{13.5\%}  & 0.0  \\
    &	    & \multicolumn{1}{c}{Average} 		       & -  & 3600.0 		        & 30.6\%  & 0.0  	   & 46.2\%  & 0.0	& 19.1\%  & 0.0 		& \textbf{18.9\%}  & 0.0     																& -  & 3600.0 	        & 30.5\%  & 0.0		  & 46.0\%  & 0.0 & 18.2\%  & 0.0 		& \textbf{17.6\%}  & 0.0  \\  
             & \multicolumn{1}{c}{84} & \multicolumn{1}{c}{Best}  & -  & 3600.0  		& 20.3\%  & 0.0		   & 44.0\%  & 0.0	& 19.3\%  & 0.0		& \textbf{18.4\%}  & 0.0     																& -  & 3600.0 		& 20.2\%  & 0.0	   & 44.4\%  & 0.0 & 17.8\%  & 0.0 		& \textbf{17.4\%}  & 0.0  \\
    &	    & \multicolumn{1}{c}{Average} 		       & -  & 3600.0 		        & 26.0\%  & 0.0  	   & 50.4\%  & 0.0	& 21.9\%  & 0.0 		& \textbf{21.4\%}  & 0.0     																& -  & 3600.0 	        & 25.9\%  & 0.0		   & 50.3\%  & 0.0 & 21.4\%  & 0.0 		& \textbf{20.5\%}  & 0.0  \\\hdashline  
    \multicolumn{3}{c}{Subgroup Average} & $\ast$ & 3600.0 & 28.6\% & 0.0 & 39.4\%  & 0.0	& \textbf{16.0\%} & 0.0 & \textbf{16.0\%} & 0.0 																									& $\ast$ & 3600.0 & 28.5\% & 0.0 & 39.3\%  & 0.0 & 15.2\% & 0.0 & \textbf{14.8\%} & 0.0  \\\hline
   \end{tabular}%
 }
 \end{center}
  \label{tabelaresultados2}%
\end{table}%
\newpage

Table~\ref{tabelaresultadoslb} compares the lower bounds obtained for each method. As it is possible to see, not only the hybrid approach achieves the best upper bounds, but also it systematically generates the best lower bounds, proving to be a very efficient and complete method. It is important to point out, that the Hybrid Lagrangian approach significantly improves the linear-relaxation bounds (column LR), even having polynomially solvable subproblems. This occurs because each time we improve the upper bound (makespan) we can transform the original problem, reducing the number of variables, allowing to improve the previous linear bound.

\begin{table}[htbp]
  \caption{Computational results for Lower Bounds.} 
 \begin{center}
  \begin{tabular}{c|cccc||cccc}
   \hline
    & \multicolumn{4}{c||}{Group 1 - Processing time  [1, 10]} & \multicolumn{4}{c}{Group 2 - Processing time  [10, 100]} \\
   \hline
   \multicolumn{1}{c|}{n} & \multicolumn{1}{c}{LR} &  \multicolumn{1}{c}{$LB_1$}  &  \multicolumn{1}{c}{$LB_2$} &  \multicolumn{1}{c||}{$LB_{HgR}$} & \multicolumn{1}{c}{LR} &  \multicolumn{1}{c}{$LB_1$}  &  \multicolumn{1}{c}{$LB_2$} &  \multicolumn{1}{c}{$LB_{HgR}$} \\
    \hline
    5 & 23.7	& 29.3 & 	33.6 & 	\textbf{34.7} & 232.8	& 288.9	& 337.4 & \textbf{347.0} \\
    10 & 39.4	& 60.9 &	66.2 & 	\textbf{70.0} & 379.1	& 635.9 & 667.6	& \textbf{706.9} \\
    20 & 70.7	& 119.7 &       140.2 &	\textbf{143.4} & 705.7	& 1203.0 & 1412.7 & \textbf{1442.4} \\
    40 & 132.9	& 225.5 &       293.5 & \textbf{294.9} & -  & 2253.2 & 2871.8 & \textbf{2883.2} \\
    60 & 193.9	& 338.0	&       \textbf{444.3} &	\textbf{444.3} & - & 3373.5	& \textbf{4452.3} & \textbf{4452.3} \\\hline
      \end{tabular}%
 \end{center}
  \label{tabelaresultadoslb}%
\end{table}%


\subsection{Results for the multi-dock problem}

\subsubsection{Instances Generation}

To investigate the performance of the proposed heuristics and the heuristics presented in literature, we generate artificial instances varying the 
number of jobs and machines, as performed in Chen and Song \cite{ChenSong2009}.
The heuristics are coded in C++ and solved using AMPL and CPLEX 12.4\footnote{all instances and codes are available at \url{https://github.com/GabrielaBragaFonseca/Cross-docking-Problems}}. 
Tests were performed on a single thread in Intel Xeon-Silver 4110 (2.1GHz/8-core/85W) with 64 GB memory and Linux operational system.

The proposed heuristics H1 and H2 are compared with the best heuristics in \cite{ChenSong2009}, called JRH and JLPTH, 
and with the CDH heuristic \cite{CotaEtAll2016}. We consider values of $n_1$ equal to $20$, $30$, $40$, $50$, $60$, $70$ and $80$, 
and values of $n_2$ equal to integers in $[0.8n_1 , 1.2n_1]$. We examine five groups of instances. The first three groups consider the same 
number of machines in each stage, $2, 4$, and $10$. The last two groups of instances use a random number of machines in each stage, 
selected from a  discrete uniform distribution $U(2,4)$ and $U(2,10)$, respectively. Processing times were generated with a discrete uniform distribution $U(10,100)$. 
A job $j \in J_1$, $J_1$ = \{$1, 2,...,n_1$\} has a probability of 50\% to belong to the set $S_{j}$ of each job 
$j \in J_2$, $J_2$ = \{$1, 2,...,n_2$\}. We generate 300 instances for each combination of number of jobs and number of machines, resulting
in a total of 10.500 instances.

\subsubsection{Computational Results}

In Table~\ref{resultadosheuristicas1} we report computational results obtained with heuristics JRH, JLPTH, CDH, H1 and H2 for large 
instances. We coded the JRH and JLPTH heuristics in C++ following the description found
in \cite{ChenSong2009}. For each combination of number of jobs and number of machines and for each heuristic we report best, average and worst 
percentage Loss over 300 instances, along with standard deviations. 
In order to evaluate the heuristics proposed we used `Loss' as the criterion to balance the quality of each heuristic algorithm, Loss = (makespan-lower bound)/lower bound, as defined by Chen and Song \cite{ChenSong2009}. 
Boldface indicates which algorithm obtained the best average result for each set of instances.

\begin{center}
\begin{scriptsize}
\begin{longtable}{lllllllll}
\caption[Computational results for large instances. The best, the average and the worst percentage Loss values are presented for each instance configuration.]{ \label{resultadosheuristicas1} Computational results for large instances. The best, the average and the worst percentage Loss values are presented for each instance configuration.}\\  \hline 
$m_i$ & \multicolumn{1}{c}{$n_1$} & \multicolumn{1}{c}{$n_2$} &  \multicolumn{1}{c}{} &  &  & Loss(\%) & & \\ \cline{5-9}
 &  &  &  & JRH & JLPTH &CDH & H1& H2 \\ \hline 
\endfirsthead
\multicolumn{9}{l}%
{{\bfseries \tablename\ \thetable{} (continued)}} \\
\hline
$m_i$ & \multicolumn{1}{c}{$n_1$} & \multicolumn{1}{c}{$n_2$} &  \multicolumn{1}{c}{} &  &  & Loss(\%) & & \\ \cline{5-9}
 &  &  &  & JRH & JLPTH &CDH & H1& H2 \\ \hline 
\endhead
\multicolumn{9}{r}{(Continued on next page)} \\
\endfoot
\hline
\endlastfoot
\multicolumn{1}{c}{2} & 20 & [16,24] & Best & 16.92 &  16.76 &  17.69 & {\bf9.02} &	9.76 \\ 
 &  &  & Average & 32.34 & 30.91 &  35.37 & {\bf20.93} & 22.20 \\   
 &  &  & Worst & 53.77 & 51.24 & 64.08 & 36.36 &	{\bf35.08} \\  
   &  &  & Std Dev  & 6.73 & 6.49 & 8.67 & 4.84 & {\bf4.74} \\ 
 & 30 & [24,36] & Best & 17.83 &  17.14 &  21.31 & {\bf11.62} &	12.46 \\ 
 &  &  & Average & 33.94 & 32.98 &  37.59 & {\bf24.25} & 25.88 \\  
 &  &  & Worst & 46.49 & 45.84 & 64.72 & 37.11 &	{\bf36.38} \\ 
  &  &  & Std Dev  & 5.39 & 5.18 & 6.96 & 4.35 & {\bf4.13} \\
 & 40 & [32,48] & Best & 22.80 & 22.29  & 24.60 & {\bf 15.54} & 17.96 \\ 
 &  &  & Average & 35.29 & 34.52 & 38.67 & {\bf26.75} &	28.36 \\  
 &  &  & Worst & 47.91 & 45.72  & 59.40 & {\bf37.13} & 39.12\\ 
   &  &  & Std Dev  & 4.34 & 4.25 & 6.41 & 3.80 & {\bf3.62} \\ 
 & 50 & [40,60] & Best & 25.58 & 25.18 & 23.65 & 19.44 &	{\bf19.36} \\ 
 &  &  & Average & 35.77 & 35.03  & 38.94 & {\bf 28.01} & 29.62\\  
 &  &  & Worst & 46.22 & 45.58 & 59.36 & {\bf 35.91} & 38.07\\  
   &  &  & Std Dev  & 3.96 & 3.91 & 6.28 & 3.32 & {\bf 3.29} \\ 
 & 60 & [48,72] & Best & 25.93 & 26.12  & 26.65 & {\bf 18.83} & 21.02\\ 
 &  &  & Average & 36.25 & 35.67 &  38.93 & {\bf 29.20} & 31.01\\  
 &  &  & Worst & 46.69 & 45.47  & 56.3 & {\bf 36.72} & 38.43\\ 
   &  &  & Std Dev  & 3.62 &  3.62 & 5.33 & 3.24 & {\bf3.12} \\ 
 & 70 & [56,84] & Best & 29.25 & 28.84  & 28.02 & {\bf 22.52} & 23.22\\ 
 &  &  & Average & 36.52 & 36.01 & 39.81 & {\bf 30.35} & 31.51\\  
 &  &  & Worst & 44.62 & 44.49 & 56.35 &  39.09 & {\bf39.07}\\ 
   &  &  & Std Dev  & 3.24 &  3.20  & 5.30 & 3.02 & {\bf 2.88} \\ 
 & 80 & [64,96] & Best & 28.84 & 28.07 &  28.01 & {\bf 22.46} & 25.42\\ 
 &  &  & Average & 36.52 & 36.12 & 39.58 & {\bf 30.99} & 32.18 \\  
 &  &  & Worst & 43.73 & 43.52 & 55.55 &  {\bf38.98} &  39.38\\  
   &  &  & Std Dev  & 2.97 &2.94  &5.13 & 2.89 & {\bf2.73} \\  \hline 
\multicolumn{3}{l}{Average 2 machines}  &  &35.23 & 34.46 & 38.41 & {\bf27.21} &28.68   \\ \hline 
\multicolumn{1}{c}{4} & 20 & [16,24] & Best & 20.90 & 17.36 & 17.56 & {\bf 11.57} & 12.99 \\
 &  &  & Average & 39.87 & 34.38 &  35.89 & {\bf 24.78} & 25.31 \\  
 &  &  & Worst & 61.61 & 52.48 & 59.75 &35.31 & {\bf 35.00}\\  
  &  &  & Std Dev  & 7.27 & 6.15& 7.43 & 4.65 & {\bf 4.33}\\ 
 & 30 & [24,36] & Best & 26.60 & 24.04 & 24.02 &  17.80 & {\bf17.47}\\ 
 &  &  & Average & 39.92 & 36.17 & 37.29 & {\bf 28.05} & 28.88\\
 &  &  & Worst & 51.71 & 46.49 &  59.00 & {\bf36.70} &  39.66 \\ 
   &  &  & Std Dev  &5.17  &  4.54 & 6.28 & 3.75 & {\bf 3.67} \\ 
 & 40 & [32,48] & Best & 26.79 & 24.49 & 22.00 & {\bf 19.42} &	20.38\\ 
 &  &  & Average & 39.68 &36.79 & 37.64 & {\bf 29.35} & 30.25\\ 
 &  &  & Worst & 50.72 & 46.13 & 54.92 & {\bf 38.28} & 39.28 \\  
   &  &  & Std Dev  & 4.31 & 3.87&  5.63 & 3.47 & {\bf 3.16} \\    
 & 50 & [40,60] & Best & 27.85 & 26.69 &  21.9 & {\bf 20.73} & 21.36 \\ 
 &  &  & Average & 39.27 & 36.71 &  37.66 & {\bf 30.29} & 31.25\\  
 &  &  & Worst & 49.65 & 47.95 & 57.31 & 39.26 & {\bf 38.59} \\  
   &  &  & Std Dev  & 3.89  &  3.68 & 5.43 & 3.29 & {\bf 3.06} \\
 & 60 & [48,72] & Best & 28.57 & 27.55 & 26.52 & {\bf 21.08} & 24.53\\ 
 &  &  & Average & 39.01& 36.96 & 37.7 & {\bf 31.07} & 32.05\\  
 &  &  & Worst & 49.20 & 45.64 & 51.43 & 39.82 & {\bf 38.62}\\ 
   &  &  & Std Dev  & 3.44 & 3.35& 4.68 &2.94 & {\bf 2.83} \\
 & 70 & [56,84] & Best & 30.86 & 30.08 & 26.71 & {\bf 23.55} & 24.86\\ 
 &  &  & Average & 38.75 & 36.99 & 38.21 & {\bf 31.52} & 32.49\\  
 &  &  & Worst & 47.80 & 44.42 & 52.69 & {\bf38.48} &  38.67 \\  
   &  &  & Std Dev  & 3.14 &2.90 & 4.57 & 2.82 & {\bf 2.62} \\ 
 & 80 & [64,96] & Best & 30.27 & 29.09 &  27.63 & {\bf 22.95} & 26.10 \\ 
 &  &  & Average &  38.39 & 36.94 & 38.33 & {\bf 32.17} & 33.12 \\ 
 &  &  & Worst & 44.81 & 43.51 &  50.97 & {\bf 38.69} & 39.38\\  
   &  &  & Std Dev  & 2.86  & 2.77 & 4.47 & 2.65 & {\bf2.60}\\  \hline    
\multicolumn{3}{l}{Average 4 machines}   &  &39.27 & 36.42 & 37.53 & {\bf29.60} & 30.48  \\  \hline 
\multicolumn{1}{c}{10} & 20 & [16,24] & Best & 12.54 & 11.59 & 5.69 & {\bf 0.00} &	4.10 \\  
 &  &  & Average & 35.38 & 24.56 & 28.23 & {\bf 13.50} &  16.49\\  
 &  &  & Worst & 55.41 & 45.50 & 61.03 & {\bf27.01} &	27.27 \\  
   &  &  & Std Dev  & 7.52 & 5.74& 8.44 & 4.79 &	{\bf 4.24}\\ 
 & 30 & [24,36] & Best & 26.91 & 18.44 & 18.53 & {\bf 12.24} &	13.86 \\ 
 &  &  & Average & 43.46 & 32.08 &  34.65 & {\bf24.14} & 24.71\\ 
 &  &  & Worst & 59.20 & 45.91 & 53.51 & {\bf 34.19} & 35.02\\ 
   &  &  & Std Dev  & 6.31 & 5.16 & 6.64 & 4.60 & {\bf 4.46} \\      
 & 40 & [32,48] & Best & 28.06 & 22.54 & 27.20 & {\bf 15.96}	& 18.67 \\ 
 &  &  & Average & 49.17 & 38.74 & 40.87 &   {\bf 31.79} & 32.10 \\  
 &  &  & Worst & 63.17 & 48.35 & 57.14 & 41.04 & {\bf 40.96} \\  
   &  &  & Std Dev  & 6.04  & 4.87 & 6.21 & 4.33 & {\bf 4.21}\\ 
 & 50 & [40,60] & Best & 36.50 & 28.39  & 27.36 & 26.04 & {\bf 25.41}\\ 
 &  &  & Average & 50.62 & 42.04 &  42.83 & 35.96 &{\bf  35.81}\\ 
 &  &  & Worst & 63.10 & 51.89 & 57.84 & 43.15 & {\bf 41.23}\\  
   &  &  & Std Dev  & 4.79 & 3.77 & 5.27 & 3.28 & {\bf 3.20} \\
 & 60 & [48,72] & Best & 38.27 & 32.56 & 30.42 & {\bf 27.74} & 27.91 \\
 &  &  & Average & 49.40 & 41.85 &  41.97 & 36.70 & {\bf 36.44}\\  
 &  &  & Worst & 59.34 & 48.44& 57.24 & 42.09 & {\bf 41.97} \\  
   &  &  & Std Dev  & 4.04 & 3.11 & 4.79 & 2.71 &	{\bf 2.67} \\ 
 & 70 & [56,84] & Best & 37.66 & 32.57 &  29.75 & {\bf28.49} & 29.02 \\  
 &  &  & Average & 47.62 & 41.21 & 41.24 & 36.53 & {\bf 36.24} \\  
 &  &  & Worst & 59.36 & 46.98 & 52.51 & 42.02 & {\bf 41.12}\\  
   &  &  & Std Dev  & 3.55 & 2.63 & 4.57 & 2.50 &	{\bf 2.34} \\   
 & 80 & [64,96] & Best & 37.61 & 33.20 & 30.17 &  28.94 & {\bf28.30}\\ 
 &  &  & Average & 45.92 & 40.44 & 40.82 & 36.34 & {\bf36.12}\\  
 &  &  & Worst & 53.71 & 45.87 & 51.81 & 41.30 & {\bf 41.29}\\  
   &  &  & Std Dev &3.11 & 2.50 & 4.38 & 2.46& {\bf 2.32} \\ \hline    
 \multicolumn{3}{l}{Average 10 machines} &  &45.94 & 37.28 & 38.66 & {\bf30.71} & 31.13\\ \hline  
 \multicolumn{1}{c}{$U(2,4)$} & 20 & [16,24] & Best & 4.57 & 4.25 & 9.02 & 3.85 & {\bf2.32}\\  
 &  &  & Average & 31.33 & 27.94 &  36.48 & {\bf 18.62} & 19.36\\  
 &  &  & Worst & 59.86 & 46.01 & 86.91 & {\bf 33.61} & 35.40\\  
   &  &  & Std Dev  &9.89  & 8.82 & 14.55 & 6.75 & {\bf 6.69}\\ 
 & 30 & [24,36] & Best & 9.57 & 6.79 & 10.14 & {\bf4.82} &	5.28 \\
 &  &  & Average & 32.07 & 29.70 & 38.14 & {\bf 21.86} & 22.68\\  
 &  &  & Worst & 51.89 & 47.81 & 85.47 & 38.32  & {\bf 36.28}\\  
   &  &  & Std Dev  & 8.26 & 7.63  & 13.76 & 6.60 & {\bf6.40}\\ 
 & 40 & [32,48] & Best & 10.90 & 9.54 & 10.61 & {\bf 3.00} &	6.57\\ 
 &  &  & Average & 32.33 & 30.58 & 38.18 & {\bf 23.72} & 24.67\\  
 &  &  & Worst & 47.42 & 44.74 &  79.10 & 37.64 & {\bf 36.84}\\  
   &  &  & Std Dev  &7.81  & 7.34  & 13.05 & {\bf6.52} & 6.57\\ 
 & 50 & [40,60] & Best & 13.26 & 12.21 & 12.64 &  {\bf8.33} & 10.21\\ 
 &  &  & Average & 32.65 & 31.27  & 38.69 & {\bf 25.04} & 26.04\\  
 &  &  & Worst & 48.53 & 45.11  & 76.62 & 38.96 & {\bf 37.57}\\  
   &  &  & Std Dev  & 7.70 & 7.34 & 12.63 & {\bf 6.53 }& 6.75\\ 
 & 60 & [48,72] & Best & 11.15 & 11.04  & 11.83 & {\bf 8.34} &	9.44 \\
 &  &  & Average & 32.74 & 31.51 & 38.74 & {\bf 25.72} &	26.84\\  
 &  &  & Worst & 47.14 & 43.42 & 81.55 & {\bf 36.36} & 38.58\\  
   &  &  & Std Dev  & 6.97 & 6.75 & 12.09 & {\bf6.17} &	6.28\\ 
 & 70 & [56,84] & Best & 16.49 & 14.56 & 16.39 & {\bf9.94} &	11.67 \\ 
 &  &  & Average & 32.36 & 31.40  & 39.3 & {\bf 26.16} &	27.32\\  
 &  &  & Worst & 45.60 & 43.95 & 77.17 & {\bf 37.74} & 39.01\\  
   &  &  & Std Dev  & 6.81  & 6.58 & 11.96 & {\bf 6.01} & 6.19 \\ 
 & 80 & [64,96] & Best & 15.01 & 14.13 & 15.48 & {\bf 11.03} & 13.05 \\ 
 &  &  & Average & 32.53 & 31.67 & 39.57 & {\bf 26.95} &	27.71\\  
 &  &  & Worst & 43.92 & 43.37 & 77.62 & {\bf 37.52} & 37.56\\  
   &  &  & Std Dev  & 6.57 & 6.41 & 12.05 & {\bf5.94} & 6.02\\  \hline    
\multicolumn{3}{l}{Average $U(2,4)$ machines} &  & 32.29 &  30.58 & 38.44 & {\bf 24.01} & 24.94\\ \hline
  \multicolumn{1}{c}{$U(2,10)$} & 20 & [16,24] & Best & 3.40 & 3.40 & 2.11 & {\bf 0.15} & 0.46 \\ 
 &  &  & Average & 27.78 & 21.25 & 35.94 & {\bf 13.15} & 13.77\\  
 &  &  & Worst & 55.42 & 43.67  & 103.55 & 34.38 & {\bf 31.48}\\  
   &  &  & Std Dev  & 11.45 & 9.18 & 22.69 & 7.57 & {\bf 7.24} \\ 
 & 30 & [24,36] & Best & 4.35 & 1.97 & 4.09 & {\bf 1.18} & {\bf 1.18} \\ 
 &  &  & Average & 31.67 & 26.05 & 40.02 & 19.44 & {\bf 19.43}\\ 
 &  &  & Worst & 61.61 & 45.32 & 97.89 & 38.78 & {\bf37.32}\\  
   &  &  & Std Dev  & 12.13  & 10.24  & 21.41 & 8.95 & {\bf 8.75}\\ 
 & 40 & [32,48] & Best & 7.67 & 4.04 & 4.50 & {\bf 2.14} & 2.72\\ 
 &  &  & Average & 33.90 & 29.08 & 40.85 &  23.07 & {\bf22.96}\\ 
 &  &  & Worst & 60.57 & 48.93 & 96.09 & 40.26 & {\bf 39.74}\\ 
   &  &  & Std Dev  & 13.05 &  10.85 & 19.21 & 9.63 & {\bf 9.58}\\   
 & 50 & [40,60] & Best & 8.09 & 6.88 & 4.61 & {\bf 3.58} & 3.77\\ 
 &  &  & Average & 33.64 & 29.66 & 41.28 & {\bf 24.30} & 24.33\\  
 &  &  & Worst & 60.51 & 50.56 & 100.21 & {\bf 40.83} & 40.86\\ 
   &  &  & Std Dev  & 12.11 & 10.53  & 18.97 & 9.46 & {\bf 9.41}\\    
 & 60 & [48,72] & Best & 8.81 & 6.33 &  6.36 &  4.54 &	{\bf4.45} \\ 
 &  &  & Average & 20.78 & 29.63 & 41.03 & {\bf 24.64} & 24.87\\  
 &  &  & Worst & 55.67 & 47.11 & 98.85 & {\bf 40.11} & 41.00\\ 
   &  &  & Std Dev  & 11.55 & 10.18 & 18.76 & {\bf 9.38} & 9.39 \\
 & 70 & [56,84] & Best & 7.49& 6.18 &  6.63 & {\bf 4.98} & 5.18\\ 
 &  &  & Average & 32.01 & 29.25 & 40.85 & {\bf 24.92} & 25.08\\ 
 &  &  & Worst & 51.53 & 47.36 & 93.43 &  40.86 & {\bf40.08}\\  
   &  &  & Std Dev  & 11.07 & 9.87 &18.56 & {\bf 9.28} & {\bf 9.28} \\ 
 & 80 & [64,96] & Best & 8.26 & 7.51 & 5.94 & {\bf 5.05} & 6.07\\ 
 &  &  & Average & 31.49 & 29.10 & 41.10 &{\bf  25.16} & 25.32\\ 
 &  &  & Worst & 50.58 & 45.83 & 92.11 & 41.67 & {\bf 40.27}\\ 
   &  &  & Std Dev  &  10.57&	9.63&	18.53&	{\bf9.01}& 9.10\\ \hline 
\multicolumn{3}{l}{Average $U(2,10)$ machines} &  & 31.93& 27.72& 40.15 & {\bf 22.10} & 22.25\\ \hline 
Total average &  &  &  & 36.93 & 33.29 & 38.64 & {\bf26.73} & 27.50\\ \hline
\end{longtable}
\end{scriptsize}
\end{center}

Results show that for all instances sizes, the proposed heuristics H1 and H2 performs better than the heuristics JRH, JLPTH and CDH, for all tested cases (Best, Average, Worst and Std Dev). 
It is noted that the total average Loss value is considerably lower for H1 and H2. 

{In this experiment, we are comparing two types of methodologies, on the one hand, we have simple polynomial time heuristics, they are very easy to code, and they are very fast. On the other hand, we have a more complicated method, that will undoubtedly need specific expertise to be implemented.  Still, this methodology is very fast, and able to present an average of 10\% of makespan benefit independently of the instance tested, reaching in some cases near 20\%. It is important to quantify that in an annual operation of 250 days, a 10\% difference in the makespan could represent a benefit of 25 days per year.}

The heuristic H1 presents a better average result than H2. Furthermore, when the number of machines increases (easier instances) the H2 heuristic performs slightly better than H1 (10 machines). However, when the number of machines decreases (more difficult instances) the performance of the H1 becomes better. 

\begin{figure}[h]
\centering
\includegraphics[scale=0.5]{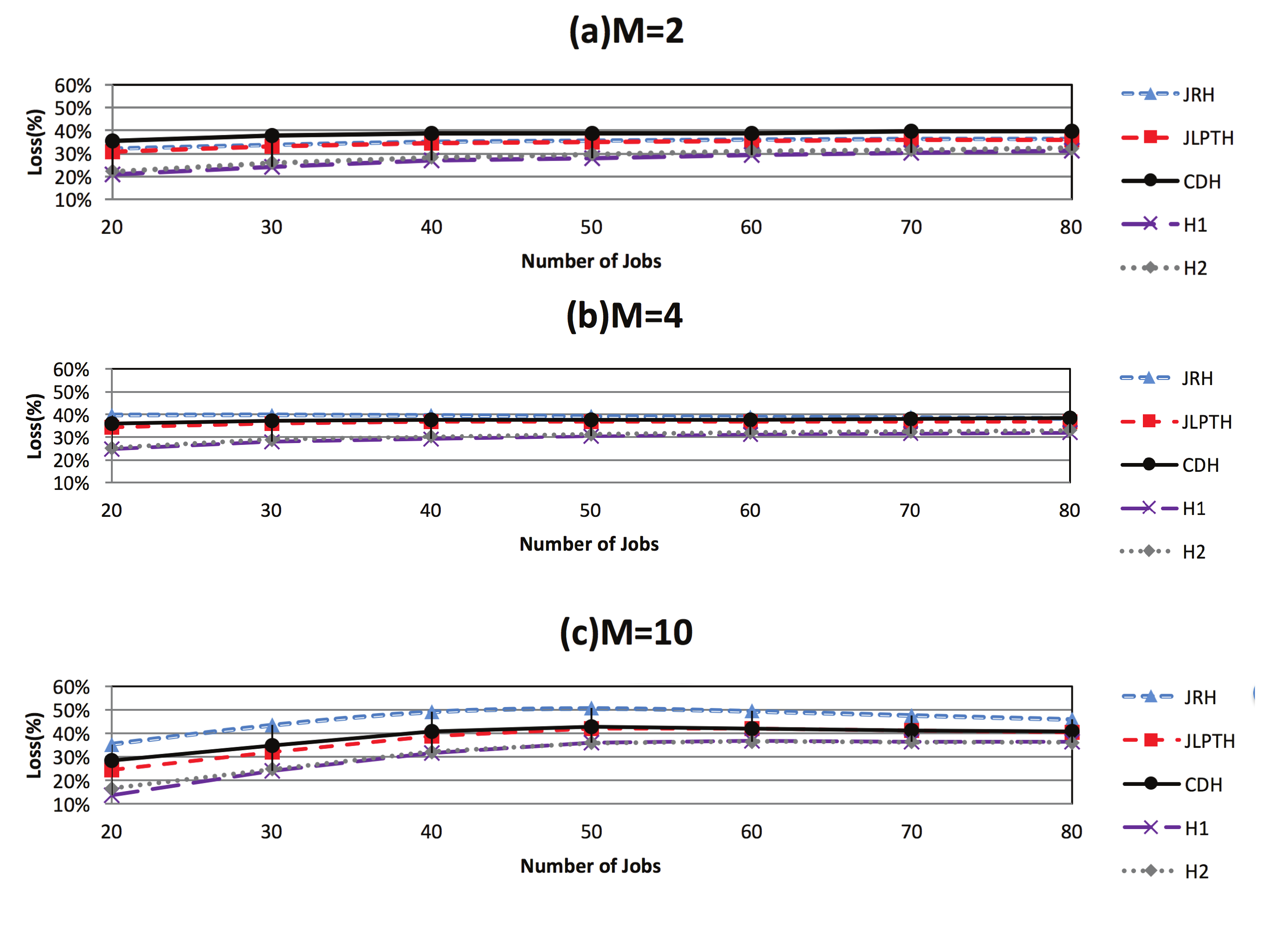} \vspace{-1.0 cm}
\caption[{Average Loss value of each heuristic considering the same number of machines at each stage}]
{{ Average Loss value of each heuristic considering the same number of machines at each stage: (a) $M=2$, (b)$M=4$, (c) $M=10$.}}
\label{graficos1}
\end{figure}
\begin{figure}[h]
\centering
\includegraphics[scale=0.5]{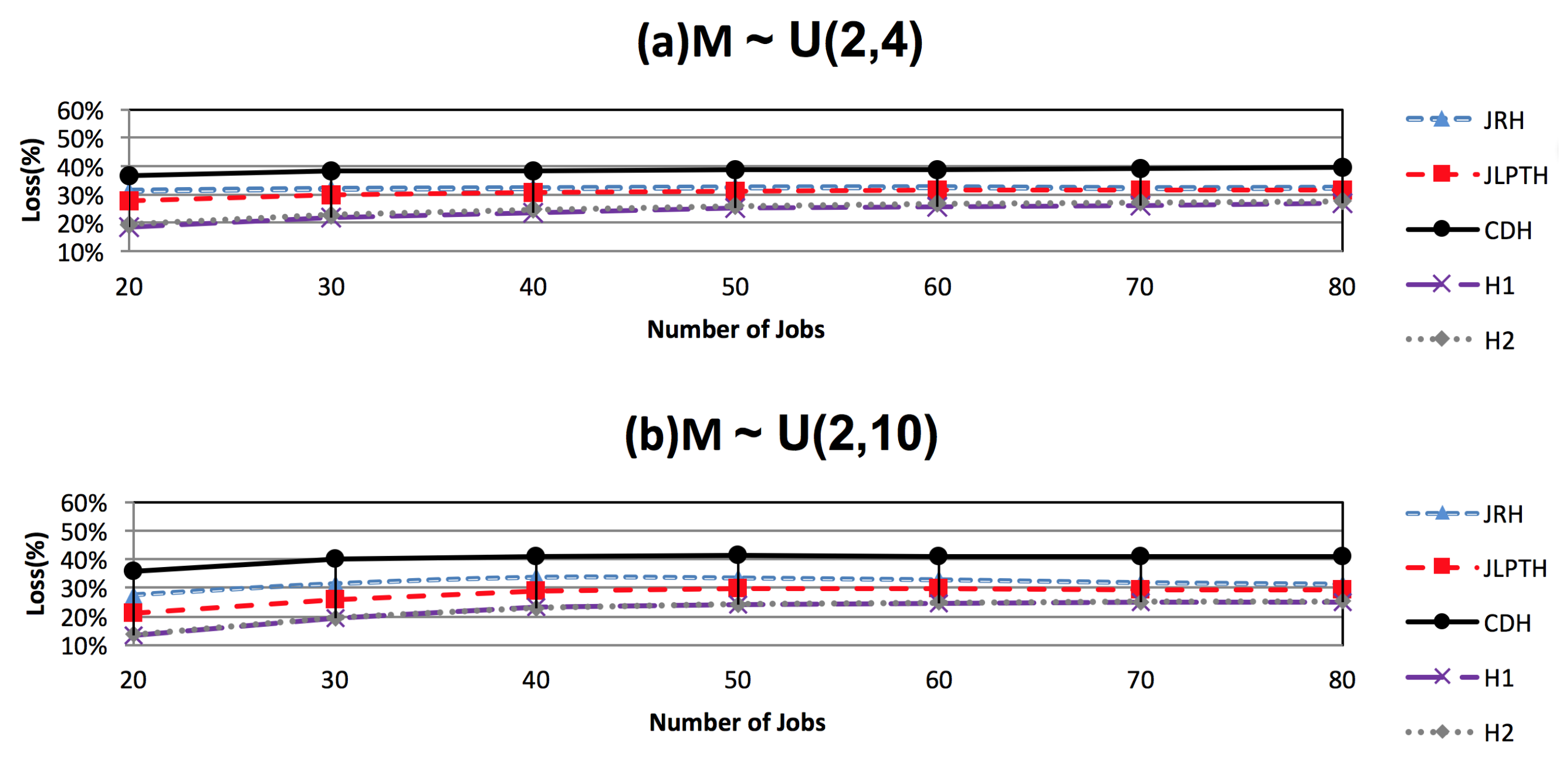} \vspace{-1.0 cm}
\caption[{Average Loss value of each heuristic considering a random number of machines at each stage} ]
{{Average Loss value of each heuristic considering a random number of machines at each stage: (a) $M \sim U(2,4)$, (b) $M \sim U(2,10)$.}}
\label{graficos2}
\end{figure}

Figures \ref{graficos1} and \ref{graficos2} summarizes the results. As previously mentioned, we can see that H1 and H2 present the lowest values in all cases. In the Figure, it is possible to see that H1 dominates the results having better Loss values in most cases. Another interesting point is that the proposed heuristics H1 and H2 have similar behavior to JRH heuristic, this may be because these heuristics arrange jobs in the second stage in increasing order of ready times.

\newpage

\newpage

\section{Conclusions}

In this work, we develop efficient ways to solve the cross-docking flow shop scheduling problem. We analyze a time-indexed formulation and a Hybrid Lagrangian Metaheuristic Framework.  As expected, the performance of the MIP is strongly dependent on the instances size, being not able to solve the problem with medium and large dimensions. 

For the 2-dock case, even having polynomially solvable Lagrangian subproblems, through a series of cuts, the method improves the model linear-relaxation bound. The Hybrid Lagrangian Metaheuristic shows an efficient performance, obtaining tight bounds in reduced computational time.

The heuristics proved to work very well with the Lagrangian approach, finding good solutions, and outperforming previous results. The heuristic H2 gives the best average GAP and Loss results followed by the heuristic H1. 

In the generalized version, the subproblems of the Lagrangian relaxation are NP-hard, and we work with their linear relaxations. The Lagrangian heuristics obtain excellent results improving the performance of previous efforts. As future research directions, we are working on new strategies for solving their subproblems and how to integrate them into the Lagrangian framework.

\section*{Acknowledgments}
This research was partially funded by FAPEMIG, CAPES and CNPq, Brazil. MGR acknowledges support from FUNDEP.

\section*{References}
\bibliography{ref-1}

\begin{appendices}

\newpage
\section{Considerations about WSPT-TRD}\label{proof1}

\begin{figure}[httb] 
    \centering
    \fbox{\includegraphics[scale=0.3]{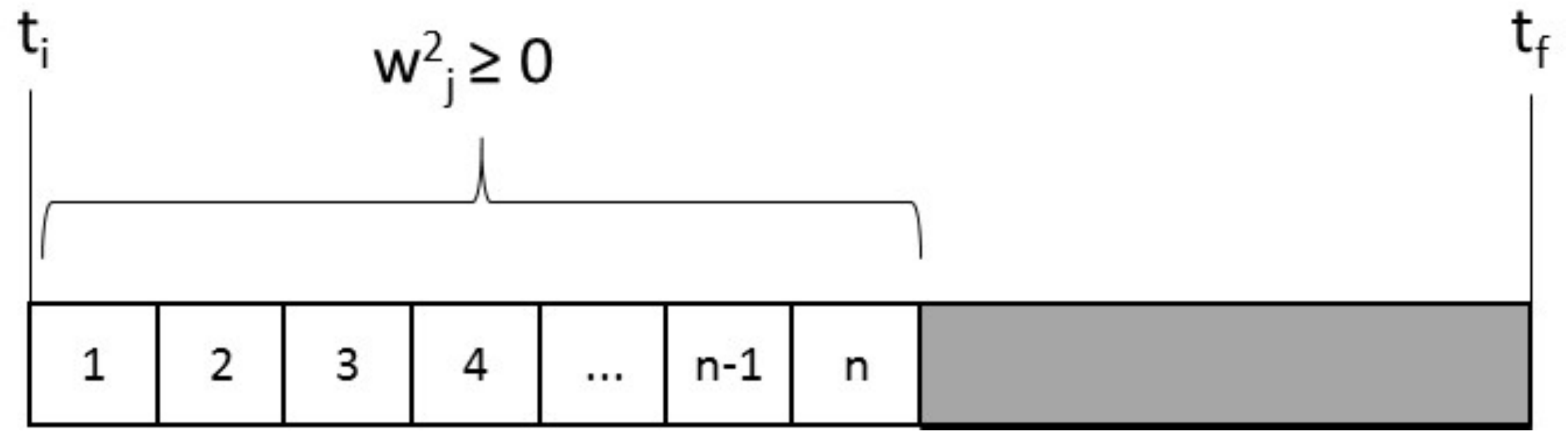}}\\
    \caption{Sequence with positive weights}
    \label{F1Proof}
\end{figure}

\begin{figure}[httb] 
    \centering
    \subfloat[Initial Schedule]{\label{F2Proof}\fbox{\includegraphics[scale=0.3]{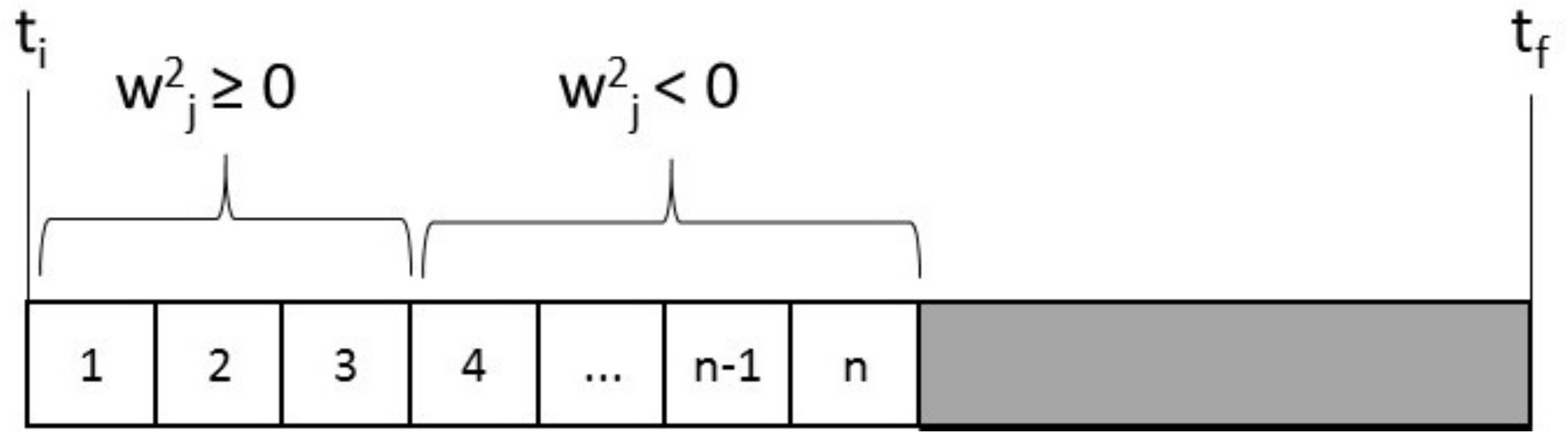}}}\\
    \subfloat[Optimal Schedule]{\label{F3Proof}\fbox{\includegraphics[scale=0.3]{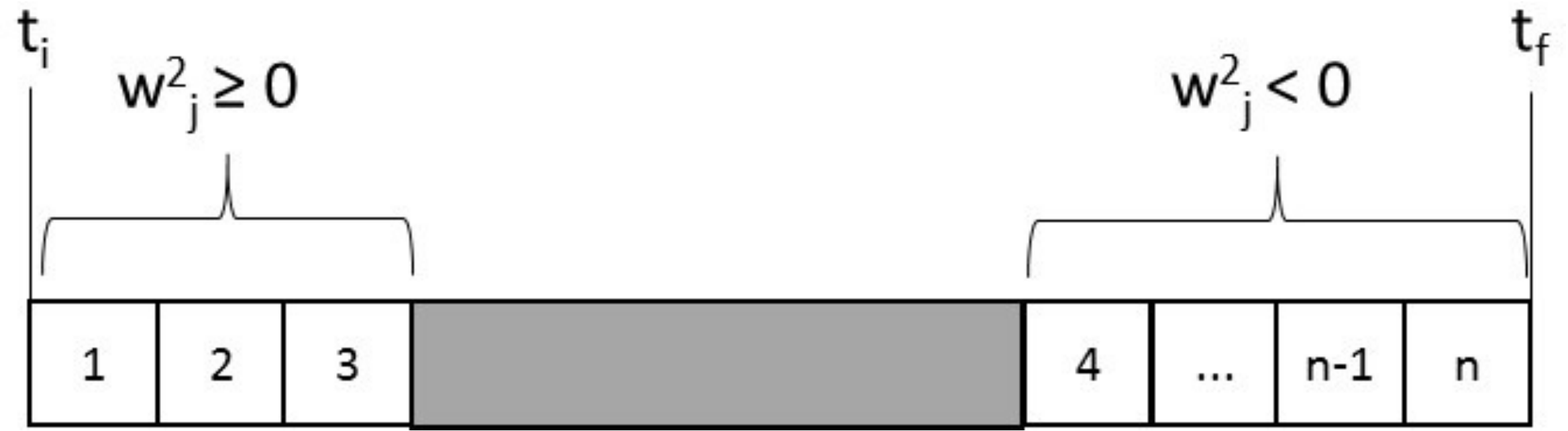}}}
    \caption{Sequence with negative and positive weights}
    \label{FigBest}
\end{figure}

\begin{thm}
Considering the problem $1||{C_{max}} +\sum I_j W_j$, the algorithm WSPT-TRD obtains an optimal solution.
\end{thm}
\begin{prf}

\noindent 
The objective function has two criteria, $C_{max}$ and $\sum I_j W_j$.
Regardless of the situation, $C_{max}$, will always aim to allocate the jobs as soon as possible. 
However, the $\sum I_j W_j$ depend on the weights values. Therefore, the proof is divided into two parts.\\

 \noindent 
Part I: If $w_j^2 \geq 0$ (see Figure \ref{F1Proof}).

 \noindent 
In this case, the two criteria of the objective function are not in conflict. The criteria $C_{max}$ and $\sum I_j W_j$ aims to schedule the
jobs as soon as possible. As there is no idle time between jobs, for $C_{max}$ any sequence starting at beginning of time horizon ($t_i$)
is optimal. In this way, the only existing criteria becomes $\sum I_j W_j$ and its optimal WSPT rule is optimal for the problem.\\

 \noindent 
Part II: If $w_j^2 < 0$.

 \noindent 
In this case, the criteria are in conflict. $C_{max}$ aim to schedule the jobs at beginning ($t_i$) of the time horizon, while $\sum I_j W_j$
at the end ($t_f$).  Consider an optimal sequence generated by WSPT rule with one subset $S'$ of jobs with $ -1 < \sum w_j^2 < 0$.
Suppose, by contradiction, that this subset must be allocated at the end ($t_f$) of the time horizon. Initially the sequence is all allocated
at the beginning ($t_i$) of the time (see Figure \ref{F2Proof}). When the schedule $S'$ with $ -1 < \sum w_j^2 < 0$ is moved in one unit of the time horizon, its
$C_{max}$ increase one unit, but $\sum I_j W_j$ decrease $\sum w_j^2$. As $C_{max}$ - $\sum w_j^2  > 0$ the objective function increase and the optimality is contradicted. Therefore, the $\sum w_j^2$ must be less or equal to 1 for allocate the jobs at the end (see Figure \ref{F3Proof}). 
It should be noted that for the case $\sum w_j^2 = -1$, these jobs are indifferent.\qed

\end{prf}

\section{Sequential Parameter Optimization Toolbox (SPOT)}\label{SPOT}

Due to its high efficiency, the Sequential Parameter Optimization Toolbox (SPOT) is the defined method to determine the parameters vector of Volume 
Algorithm. A parameters vector comprises the values of different parameters of the Volume Algorithm, i.e. it represents a candidate solution
of a parameter setting. SPOT is an implementation of the Sequential Parameter Optimization (SPO), which is an iterative model-based method of tuning 
algorithms. The tuning process is based on the parameters' data, and their utility delivered by the performance of the volume algorithm. SPO performs a 
multi-stage procedure where the model is updated at each iteration with a set of new vectors and new predictions of utility in order to improve the algorithm’s
efficiency.

The goal of SPOT is the determination of good parameters settings for heuristic algorithms. It provides statistical tools for analyzing and understanding algorithm’s performance. SPOT
is implemented as a R package, and is available in the R archive network at \url{http://cran.r-project.org/web/packages/SPOT/index.html}. Further explanation 
about SPOT can be obtained from Bartz-Beielstein \cite{Bartz-Beielstein2010}, which provides an exemplification on how SPOT can be used for automatic 
and iterative tuning. Bartz-Beielstein and Zaefferer \cite{Bartz-Beielstein2012} also give an introductory overview about tuning with SPOT.

The key elements of the SPOT methodology are algorithm design ($D_a$) and problem design ($D_p$). The first one defines ranges of parameter values 
that influence the behavior of an algorithm, such as crossover rate. These parameters are treated as variables $p_a \in D_a$ in the tuning algorithm, 
where $p_a$ represents the vector of parameters settings. $D_p$ refers to variables related to the tuning optimization problem, e.g. the search space 
dimension.

SPOT is composed by two phases: the build of the model and its sequential improvement. Phase 1 determinates an 
initial designs population from the algorithm’s parameter space. The observed algorithm is run $k$ times for each design, where $k$ is the number
of repetitions performed for each parameter setting and is increased in each run. Phase 2 leads to the efficiency of the approach and is characterized 
by the following loop:

\begin{itemize}
\item update the model given the obtained data;
\item generate design points and predict their utility by sampling the model;
\item choose the best design vectors and run $k + 1$ times the algorithm for each of them;
\item new design points are added to the population and the loop restarts if the termination
criteria is not reached.

\end{itemize}

The loop simulates the use of different parameters settings, it is subject to interactions between parameters and random effects into the experiment. 
SPOT uses results from algorithm runs to build up a meta model to tune algorithms in a reproducible way. 

\subsection{Experimental setup}

Classical works about SPOT define three different layers to analyze parameter tuning. The first one is the \textit{application layer}, considered in 
our case as the two-machine flow shop problem with cross-docking constraints. The objective function and the problem parameters are defined at
this layer. The \textit{algorithm layer}, is related to the representation of the heuristic algorithm and its required parameters are the ones that 
determine the algorithm’s performance. And finally, the \textit{design layer}, where is the tuning method, tries to find good parameter settings for 
the algorithm layer.

In this way we face two optimization problems: problem solving and parameter tuning. The problem solving covers the application layer and the Volume 
Algorithm of the algorithm layer and aims to find an optimal solution for the problem. The parameters tuning uses a tuning method to
find the best parameter values for the Volume Algorithm based on lower bounds found. The fitness value is the quality measure of the first optimization, which depends on the 
problem instance to be solved. Utility is the quality measure for parameter tuning, which reflects the performance of the Volume Algorithm for a 
vector of parameters.

Thereby, the experimental setup consists of an application layer represented by twenty instances of the two-machine flow shop problem, showed on Table
\ref{amostra}. To define the region of interest (ROI) of the tuning algorithm, the type and the lower and upper bounds of the volume algorithm's 
parameters are summarized on Table \ref{limitesparam}.

\newpage
\begin{table}[ht]
  \centering
  \caption[Sampling of instances]{Sampling of instances} 
 \begin{center}
 \smallskip
  {\footnotesize  
  \begin{tabular}{c|c}
   \hline
   Type & Instances \\
   \hline
   Y1 & 5-3-4-1 \\
   Y2 & 5-4-4-1 \\
   Y3 & 5-5-4-1 \\
   Y4 & 5-6-4-1 \\
   Y5 & 5-7-4-1 \\
   Y6 & 10-6-9-1 \\
   Y7 & 10-8-9-1 \\
   Y8 & 10-10-9-1 \\
   Y9 & 10-12-9-1 \\
   Y10 & 10-14-9-1 \\   
   Y11 & 5-3-4-2 \\
   Y12 & 5-4-4-2 \\
   Y13 & 5-5-4-2 \\
   Y14 & 5-6-4-2 \\
   Y15 & 5-7-4-2 \\
   Y16 & 10-6-9-2 \\
   Y17 & 10-8-9-2 \\
   Y18 & 10-10-9-2 \\
   Y19 & 10-12-9-2 \\
   Y20 & 10-14-9-2 \\ 
   \hline 
    \end{tabular}%
}
 \end{center}
  \label{amostra}%
\end{table}%

\begin{table}[ht]
  \centering
  \caption[Parameter bounds for tuning the Volume Algorithm]{Parameter bounds for tuning the Volume Algorithm} 
 \begin{center}
 \smallskip
  {\footnotesize  
  \begin{tabular}{c|c|c|c}
   \hline
   Parameter & Lower bound & Upper bound & Type \\
   \hline
    $\pi$ & 0.0005 & 0.0100 & FLOAT \\
    MaxWaste & 5 & 30 & INT \\
    factor & 0.1 &  1.0 & FLOAT \\
    $\alpha_{max}$ & 0.20 & 0.95 & FLOAT \\
    st & 0.6 & 1.8 & FLOAT \\ 
    $\alpha$ & 0.05 & 0.90 & FLOAT \\
    yellow & 0.20 & 0.95 & FLOAT \\
    green & 1.0 & 2.0 & FLOAT \\               
   \hline 
    \end{tabular}%
}
 \end{center}
  \label{limitesparam}%
\end{table}%

\subsection{Results of the experiment}

The SPOT algorithm, implemented in R package, is connected with the volume algorithm implemented on C++ with the aid of one callString, as presented 
below. The computer used in the tests is an Intel (R) Xeon (R) CPU X5690 @ 3.47GHz with 24 processors, 132 GB of RAM, and Ubuntu Linux operating 
system.

\begin{center}
\fbox{
\parbox{12cm}{
$
callString \leftarrow paste( $\textquotedblleft ./HgR \textquotedblright$, \pi, MaxWaste, factor, \alpha_{max}, st, \alpha, yellow, green)\\
call \longleftarrow system(callString, intern=TRUE)\\
$
\#\textit{read the results}\\
$
y = read.table($\textquotedblleft Results.txt\textquotedblright$)
$}
}
\end{center}

The Table \ref{quatromelhores} presents the four best results found by SPOT method for the volume algorithm. The first column indicates the iteration 
of SPOT, the second one presents the indicated parameters vector and the third column indicate the solution value obtained by the tested instances (lower bounds). 

The best parameters vector indicated for the volume algorithm, that generated the best lower bounds is $\pi$=0.00679, MaxWaste=24, factor=0.87753, $\alpha_{max}$=0.30337, st=1.41179, 
$\alpha$=0.08300, yellow=0.48159 and green=1.56487.

\newpage
\begin{landscape}
\begin{table}[ht]
  \caption[Four best parameters settings for the volume algorithm]{Four best parameters settings for the volume algorithm} 
 \begin{center}
  \scalebox{0.60}{
  \begin{tabular}{c|cccccccc|cccccccccccccccccccc}
   \hline
    Iteration & \multicolumn{8}{c}{Parameter} & \multicolumn{20}{c}{Solution value of the instances} \\
	     & $\pi$ & MaxWaste & factor & $\alpha_{max}$ & st & $\alpha$ & yellow & green & Y1 & Y2 & Y3 & Y4& Y5 & Y6& Y7 & Y8 & Y9 &Y10 & Y11 & Y12 & Y13 & Y14 & Y15 & Y16 & Y17 & Y18 & Y19 & Y20 \\
    \hline
     501&0.00708&28&0.22871&0.33265&1.01875&0.21709&0.80291&1.13994&24&14&30&20&30&36&46&56&60&66&155&137&191&328&414&337&537&558&629&793\\
     502&0.00527&9&0.34752&0.75030&1.71397&0.73405&0.58396&1.56980&24&15&29&21&30&37&46&56&60&66&157&154&191&341&414&337&563&558&629&793\\
     503&0.00688&20&0.62855&0.83046&1.17734&0.28130&0.60572&1.40758&25&15&30&21&30&37&46&56&60&66&118&136&191&307&414&337&531&558&629&793\\
     504&0.00679&24&0.87753&0.30337&1.41179&0.08300&0.48159&1.56487&26&15&30&20&30&32&46&56&60&66&95&136&191&307&414&337&531&558&629&793\\
     \hline 
    \end{tabular}%
}
 \end{center}
  \label{quatromelhores}%
\end{table}%
\end{landscape}
\newpage

\end{appendices}

\end{document}